\documentclass[a4paper]{amsart}

\usepackage[english]{babel}
\usepackage[utf8]{inputenc}
\usepackage{amsmath}
\usepackage{amssymb} 
\usepackage{dsfont} 
\usepackage{mathrsfs}
\usepackage{graphicx}
\usepackage[colorinlistoftodos]{todonotes}
\usepackage[arrow, matrix, curve, pic]{xy} 
\usepackage[left=3.5cm,right=3cm,top=3cm,bottom=3cm]{geometry}
\usepackage{extarrows}
\usepackage{hyperref}
\usepackage[user,titleref]{zref}
\usepackage{enumitem}
\usepackage[T1]{fontenc}
\usepackage{fancyhdr}
\usepackage[utf8]{inputenc}
\usepackage[numbers,square]{natbib}

\usepackage{chngpage} 

\theoremstyle{plain}    
\newtheorem*{thmn}{Theorem}
\newtheorem{thm}{Theorem}[section]
\newtheorem{prop}[thm]{Proposition}
\newtheorem{defn}[thm]{Definition}
\newtheorem{cor}[thm]{Corollary}
\newtheorem{lem}[thm]{Lemma}

\theoremstyle{definition}
\newtheorem*{rem}{Remark}
\newtheorem{ex}[thm]{Example}

\newtheorem*{set}{Setting}
\newtheorem{setn}{Setting}

\newcommand{\proofof}[1]{\noindent\textit{Proof (of #1).}}

\newcommand{\W}{\mathcal{W}}
\newcommand{\X}{\widetilde{X}}

\newcommand{\Xs}{X_{s}}

\newcommand{\Om}{\Omega}
\newcommand{\sub}[1]{\noindent\textbf{#1.}}
\newcommand{\reduction}[1]{\noindent\textbf{Reduction #1.}}
\newcommand{\eproof}{\hfill $\Box$\newline}
\newcommand{\F}{\varPhi}

\newcommand{\s}{\sigma}
\newcommand{\codim}{\emph{codim}}

\newcommand{\bsl}{\backslash}

\newcommand{\CC}{\mathbb{C}}

\newcommand{\ZZ}{\mathbb{Z}}
\newcommand{\NN}{\mathbb{N}}
\newcommand{\QQ}{\mathbb{Q}}
\newcommand{\PP}{\mathbb{P}}
\newcommand{\RR}{\mathbb{R}}

\newcommand{\enumeraten}{\renewcommand{\labelenumi}{\hskip 1.2em \arabic{enumi}.}}
\newcommand{\enumeratet}{\renewcommand{\labelenumi}{\hskip 1.5em(\thethm.\arabic{enumi})}}
\newcommand{\enumeratel}{\renewcommand{\labelenumi}{\hskip .5em(\thelem.\arabic{enumi})}}
\newcommand{\enumerated}{\renewcommand{\labelenumi}{\hskip .5em(\thedefn.\arabic{enumi})}}

\DeclareMathOperator{\exc}{exc}
\DeclareMathOperator{\proj}{proj}
\DeclareMathOperator{\g.c.d.}{g.c.d.}

\DeclareMathOperator{\spec}{spec}

\DeclareMathOperator{\slog}{log}
\DeclareMathOperator{\Hom}{Hom}
\DeclareMathOperator{\cone}{cone}
 
\def\hq{\hspace{-0.5mm}/\hspace{-0.14cm}/\hspace{-0.5mm}}

\title[Extension Theorems for differential forms]{Extension Theorems for differential forms on low-dimensional good quotients}

\author{Stefan Heuver}

\begin{document}

\bibliographystyle{alphadin}

\begin{abstract}
In this paper we will show that the pull-back of any regular differential form defined on the smooth locus of a good quotient of dimension three and four to any resolution yields a regular differential form.
\end{abstract}

\maketitle

\tableofcontents
\newpage

\section{Introduction and main results}

\subsection{Introduction}
In algebraic geometry we analyse differential forms to determine the geometric structure of varieties. On a smooth variety usually one considers Kähler differential forms. On a singular variety however, these differential forms might behave badly near the singular locus. A natural geometric definition of a differential form on a normal variety is a differential form that is defined on the smooth locus. From an algebraic point of view these differential forms can be considered as global section of the reflexive hull of the sheaf of Kähler differential forms and are therefore called \emph{reflexive differential forms}.

Let $X$ be a normal variety over $\CC$, $\s \in H^0\big(X,\,\Om_{X}^{[p]}\big)$ a reflexive $p$-form on $X$ ($0\leq p\leq dim(X)$) and $\eta:\X\rightarrow X$ a resolution of singularities of $X$. Since certain features of the smooth case (e.g. the Kodaira vanishing and the Serre duality) do not hold for $\s$ on $X$ in general, we would like to analyse the pull-back $\eta^*(\s)$ on $\X$ instead. To this end we need $\eta^*(\s)$ to be regular. In \citep{GK14} Graf and Kov{\'a}cs have proven that for a variety $X$ with rational singularities the pull-back $\eta^*(\s)$ extends as a regular $p$-form to all of $\X$, as long as we allow logarithmic poles along the $\eta$-exceptional locus. In general it is not known that one can avoid these logarithmic poles. However, there are some results for varieties with rational singularities in special cases. 

If $p=\dim(X)$ or $p=0$ the result follows from the definition of a normal variety with rational singularities. For $p=1$, the result was proven by Graf-Kov\'acs \citep{GK14} for du Bois singularities. For $p=2$, there is a result by Namikawa \citep{Nam01}, who works on varieties with rational $\QQ$-Gorenstein singularities.  For arbitrary values of $p$, the result has been proven by van Straaten-Steenbrink \citep{vSS85} for varieties with isolated singularities and by Greb-Kebekus-Kov{\'a}cs-Peternell \citep{GKKP11} for klt-pairs. 

\begin{rem}
In all these papers the authors show that $\eta^*(\s)$ extends to the exceptional locus without poles. Therefore we call a result like this \emph{Extension Theorem}. 
\end{rem}

\noindent\textbf{Notation.}
If we can show, that for all open subsets $U\subset X$ and any reflexive $p$-form $\s \in H^0\big(U,\,\Om_{X}^{[p]}\big)$, for $0\leq p\leq \dim(X)$, the pull-back extends to $\eta^{-1}(U)$, then we say that \emph{the Extension Theorem is true for $p$-forms on $X$}. If it is true for all values of $p$ we say, that  \emph{the Extension Theorem is true for $X$}.

\begin{rem}
As \citep{GKKP11} show in Observation $1.3$ our formulation of the Extension Theorem is equivalent to the claim that $\eta_*(\Om^p_{\X})$ is a reflexive sheaf. Since reflexivity of a sheaf is a local property, this also shows that proving the Extension Theorem is a local problem.
\end{rem}

Based on the results presented above it is natural to assume that an Extension Theorem might be proven for varieties with rational singularities in general. As a matter of fact, to this day no counter example has been presented. Besides klt-pairs, who play an important role in the minimal model program, good quotients form another useful class of varieties with rational singularities. 

Let $G$ be a reductive group and $V$ a smooth irreducible $G$ variety over $\CC$. Then we can define a good quotient $\pi:V\rightarrow X:=V\hq G$. By Boutot we know that $X$ is a normal variety with rational singularities. Examples for good quotients are Geometric Invariant Theory (GIT-) quotients (see \cite{MFK94}) that are a useful tool in the study of moduli spaces. These varieties come with an extensive amount of extra structure and additional properties. In this paper we will take advantage of these extra information about the good quotient $X$ to show that for dimension $\dim(X)\leq 4$ the Extension Theorem is true for $X$. 

\begin{rem}
Note, that every open subset of a good quotient is a good quotient. To simplify the Notation presented above, we will therefore assume without loss of generality that $U=X$.
\end{rem}

\subsection{Main results} In all the Extension Theorems of this chapter $X:=V\hq G$ is a good quotient, where $G$ is a reductive group and $V$ is a smooth $G$-variety. A precise definition and the main properties of a good quotient can be found in Chapter \ref{ch.GIT}. The main result of our paper essentially says that a reflexive $p$-form on a good quotient of dimension less than or equal to $4$ lifts to a $p$-form on any log resolution.

\begin{thm}[Extension Theorem for good quotients of dimension 4 or lower]\label{t:1}
Let $X:=V\hq G$ be a good quotient of dimension $dim(X)=n\leq 4$. Then for all values $0\leq p\leq n$ we have the following result: Let $\eta:\X\rightarrow X$ be a log resolution and let $\s\in H^0\big(X,\,\Om_{X}^{[p]}\big)$ be a reflexive $p$-form on $X$. Then 
$$
\eta^*(\s)\in H^0\big(\X,\,\Om_{\X}^{p}\big).
$$

\end{thm}

\begin{rem}
The notions of reflexive differential forms and log resolutions will be discussed in Chapter \ref{ch.facts}. We use a log resolution since it provides additional information about the exceptional locus. However, as we will see in Chapter \ref{c:ind.of choice}, the result is independent of the choice of resolution.
\end{rem}

\begin{rem}
The result is only new in dimension $3$ and $4$. Since the quotient is normal, it is smooth if the dimension $dim(X)\leq 1$. Then $\X=X$ and the resolution $\eta$ is just the identity. In dimension $2$ the Theorem follows from the fact that $X$ has finite quotient singularities (see \citep[Cor.~1]{Gur}). More information on finite quotient singularities and the Extension Theorem in this case can be found in Chapter \ref{ch.easyExtThm}.  
\end{rem}

We will prove Theorem \ref{t:1} by analysing the cases $p\in\{1,2\}$ and $p=(n-1)$ separately. It is a corollary of the following two more general results: 

\begin{thm}[Extension of $1$- and $2$-forms on good quotients]\label{t:2}
Let $X:=V\hq G$ be a good quotient of dimension $dim(X)=n$. Let $\eta:\X\rightarrow X$ be a log resolution and $\s\in H^0\big(X,\,\Om_{X}^{[p]}\big)$ a reflexive $p$-form on $X$, for $p\in\{1,2\}$. Then 
$$
\eta^*(\s)\in H^0\big(\X,\,\Om_{\X}^{p}\big).
$$

\end{thm}

Theorem \ref{t:2} directly implies Theorem \ref{t:1} in dimension $3$. Furthermore, it ensures that in dimension $4$ we only have to prove the Extension Theorem for $p=3$. This case is covered by the following theorem.

\begin{thm}[Extension of $(n-1)$-forms on good quotients]\label{t:3}
Let $X:=V\hq G$ be a good quotient of dimension $n\geq 1$. Let $\eta:\X\rightarrow X$ be a log resolution and $\s\in H^0\big(X,\,\Om_{X}^{[n-1]}\big)$ a reflexive $(n-1)$-form on $X$. Assume that the Extension Theorem holds for all reflexive $p$-forms on any good quotient of dimension less than $n$. Then 
$$
\eta^*(\s)\in H^0\big(\X,\,\Om_{\X}^{n-1}\big).
$$
\end{thm}

Theorem \ref{t:1} to \ref{t:3} hold for good quotients in general and particularly when $X$ is a GIT-quotient. Since the assertions of the three Extension Theorems are local on $X$, in Chapter \ref{ch.red. ext} we will show that for the proofs of these Extension Theorems we can reduce to the case where $V$ is an affine $G$-variety (or $V$ is a vector space with a linear action of $G$) and $X$ is the induced affine GIT-quotient. These reductions are crucial for the success of this paper since they allow us to use properties of certain GIT-quotients such as Luna's Slice Theorem (see Chapter \ref{ch. GITred}), the existence of the partial resolution of Kirwan (see Chapter \ref{ch. Kirwan}) and also provide information about the types of singularities (see Chapter \ref{ch:singcodim2}) of a good quotient. As a consequence, the proofs of the Theorems \ref{t:1} to \ref{t:3} only work for our precise set-up and cannot be generalised for arbitrary varieties with rational singularities. Nevertheless, using this special properties of GIT-quotients also means that the proofs presented in this paper are much more concrete and can provide additional insight.

\subsection{Outline of the paper} 

There are mainly two approaches to study the extension of reflexive differential forms to a resolution. The first is by studying Hodge-theoretic methods, the other is by using the residue sequence. In this paper we have to use both methods to prove Theorem \ref{t:1}. The paper is divided into three parts. 

In Part I we introduce notation and recall the most important definitions and standard facts. We will present important reductions, that are useful to analyse good quotients and to prove the Theorems \ref{t:1} to \ref{t:3}. We then prepare the proof of Theorem \ref{t:2} by introducing Hodge-theoretic methods and show that cutting down an affine GIT-quotient with hyperplanes results in a variety that is still an affine GIT-quotient.

In Part II we use a result by Gurjar \citep[Cor.~1]{Gur} and Part I to show that a good quotient has finite quotient singularities in codimension $2$. We then prove Theorem \ref{t:2} by modifying a result of Namikawa \citep[Prop.~3]{Nam01}. 

In Part III we first introduce the partial resolution algorithm of Kirwan and show that by Luna-Richardson \citep{LR79} an affine GIT-quotient $X$ of a vector space by a group (that acts linearly on the vector space) has a representation $X:=V\hq G$ such that the algorithm of Kirwan can be used. Then we prove the existence of a residue sequence on certain pairs of good-quotients and show that the assertion of Theorem \ref{t:3} holds to be true in every step of the partial resolution of Kirwan. Finally, we prove Theorem \ref{t:3}.

\vspace{7pt}
\noindent\textbf{Acknowledgements.} The author would like to thank his advisor Prof. Dr. Daniel Greb for introducing him to the field of good quotients and the Extension Theorem and for the weekly discussions on the topic. The Author would also like to thank Prof. Dr. Stefan Kebekus for a motivating discussion in Freiburg that resulted in the final version of the proof.  

\part{Preliminaries}

In this paper a variety is an integral separated scheme of finite type over $\CC$. Varieties are also assumed to be irreducible. A vector space is an affine variety isomorphic to $\CC^n$ for some $n\in \NN$. In particular a vector space is always assumed to be finite dimensional.

\section{Differential forms, resolutions and singularities}\label{ch.facts} 

\begin{defn}[Reflexive (logarithmic) differential forms {\citep[Not.~2.16]{GKKP11}}]
Let $X$ be a normal variety and $D$ a reduced divisor on $X$. For $0\leq p\leq \dim(X)$ let $\Om^p_X$ and $\Om_{X}^{p}\big(\slog D\big)$ be the sheaves of Kähler $p$-forms and logarithmic Kähler $p$-forms on $X$. Then $\Om^{[p]}_X:=(\Om^p_X)^{\vee\vee}$ and $\Om_{X}^{[p]}\big(\slog D\big):=\big(\Om_{X}^{p}\big(\slog D\big)\big)^{\vee\vee}$ are called \emph{sheaf of reflexive $p$-forms} and \emph{sheaf of reflexive logarithmic $p$-forms on $X$}, where $(~\cdot~)^{\vee\vee}$ is the double dual.
\end{defn}

\begin{rem}
For more details on reflexive sheaves the reader is referred to {\citep[Ch.~2D]{GKKP11}}. Independently, we would like to make some important remarks.
\enumeraten
\begin{enumerate}
\item Logarithmic differential forms are discussed in \citep[Ch.~11c]{Iit82}.
\item By definition the sheaves $\Om^{[p]}_X$ and $\Om_{X}^{[p]}\big(\slog D\big)$ are reflexive and in particular torsion-free.
\item Consider the pair $(X,\,D)$ and its regular part  $U:=(X,\,D)_{\mathrm{sm}}$ (see \citep[Def.~2.2,~2.4]{GKKP11}). If we denote by $i:U\hookrightarrow X$ the embedding into $X$, then $\Om_{X}^{[p]}\big(\slog D\big)=i_*\big(\Om_{U}^{p}\big(\slog D_{\mid U}\big)\big)$ (see \citep[Not.~2.16]{GKKP11}).
\end{enumerate}
\end{rem}

\begin{defn}[Reflexive relative differential forms]
Let $\varPsi:X\rightarrow T$ be a morphism from a normal variety $X$ to a smooth variety $T$. For $0\leq p\leq dim(X)$, let $\Om^p_{X/T}$ be the sheaf of relative $p$-forms on $X$ (see \cite[II.~Prop.~8.11]{Har77}). Then by $\Om^{[p]}_{X/T}=(\Om^p_{X/T})^{\vee\vee}$ we denote the \emph{sheaf of reflexive relative $p$-forms on $X$}.
\end{defn}

\begin{defn}[Rational singularities {\citep[Def.~5.8]{KM98}}]
Let $X$ be a normal variety. We say that $X$ has \emph{rational singularities} if there exists a proper birational map $f:Y\rightarrow X$ from a smooth variety $Y$, such that $R^i f_*\mathcal O_Y=0$ for $i>0$. 
\end{defn}

\begin{ex}
The typical example is the singularity of the quadric cone given by the equation $x^2+y^2+z^2=0$.
\end{ex}

\begin{prop}[{\citep[Lemma~5.12]{KM98}}]\label{nforms}
Let $X$ be a variety of dimension $n$ and $f:Y\rightarrow X$ a resolution of singularities of $X$. Then $X$ has rational singularities if and only if $X$ is Cohen-Macaulay and $f_*\Om^{n}_Y=\Om^{[n]}_X$.
\end{prop}

\proof
See \citep[Lemma~5.12]{KM98}.
\eproof

\begin{rem}[Extension of $0$ and $n$-forms]
The equation $f_*\Om^{n}_Y=\Om^{[n]}_X$ implies that every reflexive $n$-form on $X$ pulls back to all of $Y$ without poles. Since normality of $X$ implies that $f_*\mathcal O_Y=\mathcal O_X$, the same is true for $0$-forms on $X$.
\end{rem}

In many cases of practical interest the traditional definition of a resolution (see \cite[Not.~04.(9)]{KM98}) is not enough, since we would like to have addition information on the exceptional locus.

\begin{defn}[Log resolution {\citep[Def.~2.12]{GKKP11}}]
Let $X$ be a variety. A \emph{log resolution} is a surjective birational morphism $\eta:\X\rightarrow X$ such that
\enumerated
\begin{enumerate}
  \item the variety $\X$ is smooth,
  \item the $\eta$-exceptional set $\exc(\eta)$ is a divisor with simple normal crossing (snc).
\end{enumerate}
We call $\eta:\X\rightarrow X$ a \emph{strong log resolution} if the following property holds in addition 
\begin{enumerate}[resume]
  \item The morphism $\eta_{\mid \eta^{-1}(X_{\mathrm{sm}})}:\eta^{-1}(X_{\mathrm{sm}})\rightarrow X_{\mathrm{sm}}$ is an isomorphism.
\end{enumerate}
By abuse of notation we will sometimes call $\X$ resolution of $X$.
\end{defn}

\begin{rem}
By Hironaka's Theorem (c.f. \citep[pp.~3-4]{KM98}) (strong) log resolutions exist in our setting. The definition of an snc divisor can be found in \citep[0.4.(8)]{KM98}. Since the $\eta$-exceptional set is a divisor, we will sometimes call it \emph{$\eta$-exceptional divisor}.
\end{rem}

\begin{prop}[Extension with log poles {\citep[Thm.~4.1]{GK14}}]\label{logpoles}
Let $X$ be a normal $n$-dimensional variety with rational singularities and let $\eta:\X\rightarrow X$ be a log resolution with exceptional divisor $E:=\exc(\eta)$. Let $\s\in H^0\big(X,\,\Om_{X}^{[p]}\big)$ be a reflexive $p$-form on $X$, for $0\leq p\leq n$. Then 
$$
\eta^*(\s)\in H^0\big(\X,\,\Om_{\X}^{p}\big(\slog E\big)\big).
$$
We say that $\eta^*(\s)$ has \emph{log-poles along $E$}.
\end{prop}		

\proof
By Theorem $S$ in \citep{Kov99} the variety $X$ is Du Bois. Thus the result follows from \citep[Thm.~4.1]{GK14}.
\eproof

\section{Good, geometric and GIT-quotients}\label{ch.GIT}
In this Chapter we want to give a short introduction on Mumfords' Geometric Invariant Theory, \citep{MFK94}. The chapter covering the affine case is based on \citep{Kra84}. Our discussion of the projective case is based on \citep{Dol03}. After these introductions we deal with the properties of a GIT-quotient of $\CC^n$ in detail and present techniques to reduce to this case. Throughout this chapter $G$ will be a reductive group. We start with a general definition of good and geometric quotients.

\begin{defn}[Good and geometric quotient {\citep[Def.~1.5, Def.~1.6]{Ses72}, \citep[p.~143]{BS91}}]\label{d:good quot}
Let $G$ be  a reductive group and $V$ a $G$-variety. A morphism $\pi: V\rightarrow X$, where $X$ is an algebraic space, is called \emph{good quotient (of $V$ by (the action of) $G$)} if the following properties are fulfilled:
\enumerated
\begin{enumerate}
\item $\pi$ is $G$-invariant, 
\item $\pi$ is affine,
\item $\mathcal O_X\cong\pi_*(\mathcal{O}_V^G)$.
\end{enumerate}
It is called \emph{geometric quotient} if the following property holds in addition:
\begin{enumerate}[resume]
\item For each point $x\in X$, the fiber $\pi^{-1}(x)$ is a $G$-orbit.
\end{enumerate}
By abuse of notation we call $X$ good or geometric quotient and use the notation $X=V\hq G$ or $X=V/G$ respectively. 
\end{defn}

\begin{rem}[Properties of good and geometric quotients {\citep[Ch.~1]{BS91}}]\label{r:properties of good quotients} 
Recall from \citep[Ch.~1]{BS91} that a subset $W\subset V$ is called \emph{$G$-saturated (in $V$)}, if for all $v\in W$ the closures of the orbit $G(v)$ in $V$ and $W$ coincide. The good quotient $X$ has the following useful properties:
\enumeraten
\begin{enumerate}
\item Let $W\subset V$ be open and $G$-saturated in $V$. Then $\pi(W)\subset X$ is open in $X$ and $\pi_{\mid W}:W\rightarrow \pi(U)$ is a good quotient. Moreover, for any (open) sub-variety $U\subset X$ the preimage $\pi^{-1}(U)$ is a $G$-saturated (open) subset of $V$.
\item Let $W_1$, $W_2\subset V$ be two disjoint, closed and $G$-invariant subsets of $V$. Then the images $\pi(W_1)$ and $\pi(W_2)$ are disjoint in $X$. 
\end{enumerate}
A consequence of the second property is that every fiber of $\pi$ contains exactly one closed $G$-orbit. In the special case where $X$ is a geometric quotient, this $G$-orbit coincides with the fiber of $\pi$.
\end{rem}

\subsection{Affine GIT-quotients}
When $V$ is an affine $G$-variety there exists a intuitively accessible description of an affine GIT-quotient. In most of the proofs we will reduce to this case. 
\begin{defn}[Affine GIT-quotient {\citep[II.3.2]{Kra84}}]\label{d:affineGIT}
Let $G$ be a reductive group and $V$ an affine $G$-variety. Then the ring of $G$-invariant polynomials $\CC[V]^G$ is a finitely generated algebra. We call $V\hq G:=\spec(\CC[V]^G)$ the \emph{affine GIT-quotient of $V$ by $G$}. The natural projection $\pi:V\rightarrow V\hq G$ is called quotient map.
\end{defn}

\begin{rem}
The quotient $V\hq G$ might not be a geometric quotient, but is a good quotient (the necessary properties are discussed in {\citep[II.3.2]{Kra84}}). In particular, the quotient map is affine and surjective. On the other hand by the first property presented in the remark to Definition \ref{d:good quot}, every good quotient can locally be described as an affine GIT-quotient.
\end{rem}

\subsection{Semi-stable and stable points and GIT-quotients of projective varieties}
When $V$ is a projective $G$-variety we need to restrict ourselves to the open subsets of (semi-)stable points $V^{s}\subset V^{ss}\subset V$ to be able to define a GIT-quotient. Let us recall the definition of a (semi-)stable point.

\begin{defn}[$G$-linearisation of a line bundle {\citep[Ch.~7.1]{Dol03}}]
Let $V$ be a $G$-variety (with a $G$-action $\rho: G\times V\rightarrow V$) and let $\mathcal L$ be a line bundle on $V$. A \emph{$G$-linearisation of $\mathcal L$} is a $G$-action $\overline{\rho}:G\times \mathcal L\rightarrow \mathcal L$ such that
\enumerated
\begin{enumerate}
   \item the diagram 
   $$
	\xymatrix{
     G\times \mathcal L\ar[r]^{\overline{\rho}}\ar[d]& \mathcal L\ar[d]\\
     G\times V         \ar[r]^{\rho}       & V
    }
   $$
   commutes and
   \item the zero section of $\mathcal L$ is $G$-invariant.
\end{enumerate}
By abuse of notation we will refer to $\mathcal L$ as \emph{a $G$-linearisation on $V$}.
\end{defn}

\begin{rem}[The induced action on global sections {\citep[Ch.~7.3]{Dol03}}]
The $G$-linearisation $\overline{\rho}$ of $\mathcal L$ induces a $G$-action $\overline{\rho}'$ on the space of global sections $H^0(V,\,\mathcal L)$ of $\mathcal L$. This action is given by $$\overline{\rho}'(g,s)(v):=\overline{\rho}(g,s(\rho(g^{-1},v)))$$
for all $g\in G$, $s\in H^0(V,\,\mathcal L)$ and $v\in V$.
\end{rem}

\begin{ex}\label{e:linearisation}
Let $V$ be an affine $G$-variety and $\mathcal L:=V\times \CC$ the trivial line bundle. Then a $G$-linearisation $\overline{\rho}$ of $\mathcal L$ is given by $\overline{\rho}(g,(v,z)):=(\rho(g,v), z)$, where $g\in G$, $(v,z)\in V\times \CC$ and $\rho:G\times V\rightarrow V$ is the $G$-action on $V$. We will refer to this as the \emph{trivial $G$-linearisation on $V$}.
\end{ex}

\begin{defn}[(Semi-)stable point {\citep[Ch.~8.1]{Dol03}}]
Let $V$ be a $G$-variety and $\mathcal L$ a $G$-linearised line bundle on $V$.
\enumerated
\begin{enumerate}
\item A point $v\in V$ is called \emph{semi-stable (with respect to $\mathcal L$)} if there exists an integer $m\geq 1$ and a $G$-invariant section $s\in H^0\big(V,\,\mathcal{L}^{\otimes m}\big)^G$, such that $U_s:=\{u\in V\mid s(u)\neq 0\}\subset V$ is affine and $v\in U_s$.

\item A semi-stable point $v\in V^{ss}$ is called \emph{stable (with respect to $\mathcal L$)} if all orbits of $G$ in $U_s$ are closed and the stabiliser $G_v$ is finite.
\end{enumerate}
\end{defn}

\begin{rem}
Note, that the definitions of $V^s$ and $V^{ss}$ do not change, if we replace $\mathcal L$ by a positive tensor power of $\mathcal L$ (see \cite[Rem.~8.1.4]{Dol03}). 
\end{rem}

\begin{rem}
If $V$ is affine, then we can choose $\mathcal L:=V\times \CC$ and the trivial $G$-linearisation introduced in Example \ref{e:linearisation} on $V$. Every section $s\in H^0\big(V,\,\mathcal{L}\big)^G$ corresponds to a $G$-invariant polynomial on $V$ (because the $G$-action on the space of global sections corresponds to the $G$-action on the affine coordinate ring). If for every $v\in V$ we choose $s$ to be the section that corresponds to a $G$-invariant constant non-zero polynomial, we get $V^{ss}= V$. 
\end{rem}

\noindent\textbf{Construction of the GIT-quotient.}
In the setting of the previous definition one can construct the GIT-quotient $V^{ss}\hq G$ in the following way: We first cover $V^{ss}$ by finitely many $U_{s_i}$, where $s_i\in H^0\big(V,\,\mathcal{L}^{\otimes m_i}\big)^G$ for $i\in I$ ($I$ is a finite set) and $m_i\geq 1$. Then gluing the affine GIT-quotients $U_{s_i}\hspace{0.2pt}\hq G$, we get a quotient $X:=V^{ss}\hq G$. More details on this can be found in {\citep[Thm.~8.1]{Dol03}}. 

\begin{rem}
The quotient $V^{ss}\hq G$ is a good quotient with quotient map $\pi:V^{ss}\rightarrow V^{ss}\hq G$. The set $V^{s}\subset V^{ss}$ is a saturated open subset and $V^{s}/G$ is a geometric quotient (see {\citep[Thm.~8.1]{Dol03}}). We are mainly interested in the following two special cases:
\enumeraten
\begin{enumerate}
	\item If $V$ is projective and $\mathcal L$ is ample, then the quotient $X$ is projective and given by $X:=\proj(R^G)$, where $R:=\bigoplus_{n\geq 0}H^0\big(V,\mathcal L^n\big)$ (see \citep[Prop.~8.1]{Dol03}). 
	\item If $V$ is affine and $\mathcal L$ is the trivial linearisation, then there exists a section $s\in H^0\big(V,\,\mathcal{L}\big)^G$ that corresponds to a non-zero constant polynomial on $V$ and $U_s=V$. Thus the quotient $X$ is the affine GIT-quotient defined in Definition \ref{d:affineGIT}.
\end{enumerate}

\end{rem}

\begin{defn}[Projective GIT-quotient]\label{d:projGIT}
Let $G$ be a reductive group and $V$ a projective $G$-variety with a $G$-linearisation $\mathcal L$ such that we can construct a quotient $X:=V^{ss}\hq G$. If $\mathcal L$ is corresponds to a projective embedding $V\hookrightarrow \PP^n$ for some $n\in\NN$ (i.e. it is ample). Then we call $X$ \emph{projective GIT-quotient of $V$ by $G$}.
\end{defn}

\subsection{The case $V=\CC^n$ and $\CC^*$-actions}\label{subs:C^n}

Let $V$ be an affine variety with affine coordinate ring $A:=\CC[V]$. Let $\rho: \CC^*\times V\rightarrow V$ be any action of the group $\CC^*$ on $V$. Then $\rho$ induces an action $\rho': \CC^*\times A\rightarrow A$ on the coordinate ring $A$ via $\rho'(t,f)(v):=f(\rho(t^{-1},v))$ for all $v\in V$. For every integer $d\in \ZZ$ we consider the corresponding character $\chi_d:G\rightarrow \CC^*$ of the reductive group $G=\CC^*$ (see \cite[§1]{CLS11}). Then $A_d:=\{f:V\rightarrow \CC\mid \rho'(t,f)=\chi_d(t)\cdot f\}$ defines a grading $A=\bigoplus_{d\in \ZZ} A_d$ on $A$, where $A_k\cdot A_l\subset A_{k+l}$. We call the $\CC^*$-action on $V$ \emph{good} if $A_d=0$ for all $d< 0$ and $A_0=\CC$. In this case the only fix-point $v_0\in V$ of the $\CC^*$-action (corresponding to the maximal ideal $\bigoplus_{d>0} A_d$) is called \emph{vertex} (see \cite[Ch.~1]{Pin77}).

\begin{ex}
Let $V=\CC^n$ and $\rho: \CC^*\times\CC^n\rightarrow\CC^n$ a $\CC^*$-action on $\CC^n$ given by $\rho(t,(z_1,\dots, z_n)):=(t^{q_1}\cdot z_1,\dots, t^{q_n}\cdot z_n)$. Then $\rho$ is a good $\CC^*$-action if and only if $q_i> 0$ for all $i\in\{1,\dots, n\}$ and $\g.c.d.(q_1,\dots, q_n)=1$. The vertex of this good $\CC^*$-action is $0\in \CC^n$.
\end{ex}

\begin{lem}[Good $\CC^*$-action on GIT-quotients]\label{l:goodC}
Let $G$ be a reductive group acting linearly on $V=\CC^n$. We denote the GIT-quotient by $X:=V\hq G$. Then the good $\CC^*$-action on $V$, introduced in the previous example (with $q_i=1$ for all $i\in\{1,\dots, n\}$), induces a good $\CC^*$-action on $X$.
\end{lem}

\proof
Let $G\times V\rightarrow V \, , \, (g,v)\mapsto\varrho(g)v$ be the linear action of $G$ on $V$, where $\varrho:G\rightarrow GL(V)$ is a representation of $G$. As we have seen in the previous example, there exists a good $\CC^*$-action 
$$
m:\CC^*\times V\rightarrow V \, , \, (t,v)\mapsto m_t v
$$
on $V$, where $m_t\in GL(V)$ is the multiplication with $t\in \CC^*$. By definition, for each $g\in G$ and each $t\in \CC^*$ we have $\varrho(g) m_t=m_t\varrho(g)$ and the two actions commute. Therefore, $m$ induces a $\CC^*$-action on $X$. 
The coordinate ring of $X$ is finitely generated (since $G$ is reductive) and given by $\CC[X]:=\CC[V]^G=\bigoplus_{d\in\NN_0}\CC[X]_{(d)}$, where $\CC[X]_{(d)}=\CC[V]_{(d)}^G=(\CC[V]_{(d)})^G$ are the $G$-invariant polynomials on $V$ of degree $d$. We may see that $\CC[X]_{(0)}=(\CC[V]_0)^{G}=\CC$. Thus, $X\hq \CC^*=\{pt\}$ and the $\CC^*$-action has a unique fixpoint $x_0\in \pi^{-1}(X\hq \CC^*)$ given by the unique closed orbit. This point $x_0$ corresponds to the maximal ideal $\bigoplus_{d\in\NN_{>0}}\CC[X]_{(d)}$ and is the vertex of the good $\CC^*$-action on $X$. \eproof

\subsection{Properties and reductions}\label{ch. GITred}
Our main results are formulated for an arbitrary good quotient $X:=V\hq G$ of a smooth variety $V$ by a reductive group $G$ and we would like to analyse the properties of these quotients. Some of the properties (e.g. singularities) can be analysed locally. We already mentioned that every good quotient can locally be described as an affine GIT-quotient. In this chapter we will present a result of Gurjar \citep{Gur}, who uses Luna's Slice Theorem \citep{Lun73} to reduce from the affine case to the case $V=\CC^n$. We then recall that any good quotient of a smooth variety has rational singularities. 

\begin{lem}[Reductive stabiliser]\label{l:Mats}
Let $G$ be a reductive group and $V$ an affine $G$-variety. Assume that $v\in V$ has a closed orbit $G(v)$. Then the stabiliser $G_v$ is reductive.
\end{lem}

\proof
This result was proven by \citep{Mat60}.\eproof

The following formulation of Luna's Slice Theorem as well as more information about it can be found in \citep{Dre12}. There Dr\'{e}zet also defines $G$-morphisms (\citep[Def.~2.4]{Dre12}) and (strongly) \'{e}tale ($G$-)morphisms (\citep[Def.~4.1; Ch.~4.4]{Dre12}).

\begin{thm}[Luna's Slice Theorem {\citep[Thm.~5.3; Thm.~5.4]{Dre12}}]\label{t:LST}
Let $G$ be a reductive group and $V$ an affine $G$-variety. Let $v\in V$ be a point such that the orbit $G(v)$ is closed (and thus $G_v$ reductive (see Lemma \ref{l:Mats})). Then there exists a locally closed sub-variety $S\subset V$ of $V$, called \emph{slice}, such that 
\enumeratet
\begin{enumerate}
\item $S$ is affine and $v\in S$,
\item $S$ is $G_v$ invariant,
\item the image of the $G$-morphism $G\times_{G_v} S:=(G\times S)\hq G_v\rightarrow V$ induced by the $G$-action on $V$ is a saturated open subset $U\subset V$,
\item the restriction $G\times_{G_v} S\rightarrow U$ is a strongly \'{e}tale surjective $G$-morphism.
\end{enumerate}

If in addition $v\in V$ is smooth, we get an \'{e}tale $G_v$-invariant morphism $\varPhi:S\rightarrow T_vS$ to the tangent space to $S$ at $v$ such that $\varPhi(v)=0$, $T\varPhi_v=Id$ and such that
\begin{enumerate}[resume]
\item $T_vV=T_v(G_v)\oplus T_vS$,
\item the image of $\varPhi$ is a saturated open subset $U'\subset T_vS$,
\item the restriction $S\rightarrow U'$ is a strongly \'{e}tale surjective $G_v$-morphism.
\end{enumerate}
\end{thm}

\proof
A proof and more details on the theorem can be found in \citep[Ch.~5]{Dre12}
\eproof

\begin{rem}
Since the surjective $G$-morphism in \ref{t:LST}.$4$ is strongly \'{e}tale, it induces a surjective \'{e}tale morphism 
$$
G\times_{G_v} S\hq G\cong S\hq G_v\rightarrow U\hq G.
$$
Similarly, the $G_v$-morphism in \ref{t:LST}.$7$ induces a surjective \'{e}tale morphism 
$$
S\hq G_v\rightarrow U'\hq G_v.
$$
\end{rem}

Building on this, Gurjar proves the following result that will help us to analyse the singularities of good quotients in Chapter \ref{ch:singcodim2}. 

\begin{cor}[c.f. {\citep[Thm.]{Gur}}]\label{l:gurjar}
Let $G$ be a reductive group and $V$ a smooth affine $G$-variety. Let $X:=V\hq G$ be the induced quotient and $x\in X$ any point. Then there exists a reductive group $H$ acting linearly on some $\CC^n$, such that the analytic local ring of $X$ in $x$ is isomorphic to the analytic local ring of $\CC^n\hq H$ at the image of $0\in \CC^n$ in $\CC^n\hq H$. 
\end{cor}

\proof
The proof can be found in \citep[§1]{Gur}. \eproof

\begin{rem}
In the setting of the previous corollary, let $\pi:V\rightarrow X$ be the quotient map and $v\in \pi^{-1}(x)$ a point with closed orbit. We know that $v\in V$ is smooth, $H=G_v$ is the stabilizer of $v$ and $\CC^n$ is isomorphic to the tangent space $T_v S$ to the slice $S$ at $v$ (see Theorem \ref{t:LST}).
\end{rem}

\begin{prop}
Let $G$ be a reductive group and $V$ a smooth $G$-variety admitting a good quotient $X:=V\hq G$. Then $X$ is normal and has rational singularities.
\end{prop}

\proof
Since being normal and having rational singularities are local properties it is enough to prove them in the case where $V$ is an affine $G$-variety and $X$ is an affine GIT-quotient. The normality of $X$ follows from the fact that $V$ is normal (see \citep[Ch.~3.3, Satz~1]{Kra84}). The result that an affine GIT-quotient $X$ has rational singularities is a famous theorem by Boutot (\citep[Cor.]{Bou87}). A proof can also be found in \citep{Kov00}.
\eproof

\subsection{Principal points and finite stabiliser}

Let $G$ be a reductive group and $V$ a smooth affine $G$-variety with affine GIT-quotient $X=V\hq G$. If we choose the trivial $G$-linearisation on $V$, then all points $v\in V$ are semi-stable, $V^{ss}=V$. The stable points are exactly the $v\in V$ with finite stabiliser and closed orbit. In this chapter we want to present a method by Luna-Richardson to find an affine variety $F$ and a reductive group $\W$, such that $F\hq \W\cong X$ and $F^{s}\neq\emptyset$. More details on this method can be found in \citep{LR79}. In this chapter we will denote by $O(x)$ the unique closed $G$-orbit in the fiber of $x\in X$. 

\begin{defn}[Principal point {\citep[Def.~3.2]{LR79}}]
Let $G$ be a reductive group and $V$ an affine $G$-variety with affine GIT-quotient $X:=V\hq G$. A point $x\in X$ is called \emph{principal point}, if there exists an open neighbourhood $U\subset X$ of $x$, such that for all $x'\in U$ the following condition is fulfilled: If $O(x)$ is the closed orbit in the fiber of $x$ and $O(x')$ is the closed orbit in the fiber of $x'$. Then $O(x)$ and $O(x')$ are $G$-isomorphic homogeneous spaces. We denote by $X_{\mathrm{pr}}$ the space of all principal points.
\end{defn}

\begin{rem}
If $v\in O(x)$ and $v'\in O(x')$. Then $O(x)$ and $O(x')$ are $G$-isomorphic homogeneous spaces if and only if the stabilisers $G_x$ and $G_{x'}$ are conjugated subgroups of $G$ (see \citep[Rem.~3.3]{LR79}). The space $X_{\mathrm{pr}}\subset X$ is a dense open subspace of $X$ (see \citep[Lem.~3.4]{LR79}). 
\end{rem}

\begin{lem}
Let $H\subset G$ be a closed reductive subgroup of $G$ and $X^H:=\{x\in X\mid h.x=x$ for all $h\in H\}$ the space of all $H$-invariant points in $X$. Then the normaliser $N_G(H):=\{g\in G\mid gH=Hg\}$ is a reductive group with an induced action on $X^H$. 
\end{lem}

\proof
The proof of the first assertion can be found in \citep[Lem.~1.1]{LR79}. The second assertion follows directly from the definitions.
\eproof

\begin{prop}[{\citep[Thm.~4.2]{LR79}}]\label{l:lunarichardson}
Let $G$ be a reductive group and $V$ an affine $G$-variety with affine GIT-quotient $X:=V\hq G$. Let $x\in X_{\mathrm{pr}}$ be a principal point and $a\in O(x)\subset V$. Then $H:=G_a$ is reductive by Lemma \ref{l:Mats}. Set $\W:=N_G(H)/H$ and $F:=X^H$. Assume that $F\hq \W$ is irreducible. Then $X\cong F\hq \W$.
\end{prop}

\proof
The proof can be found in \citep[Thm.~4.2]{LR79}.
\eproof

\begin{rem}
If $V$ is smooth, the restriction that $F\hq \W$ is irreducible is not necessary. A detailed explanation can be found in \citep[Rem.~4.6]{LR79}. If $V$ is a vector space and $G$ acts linearly on $V$. Then one can easily check that $F$ is a vector space as well and $\W$ acts linearly on $F$. Since the vector space $F$ is irreducible, $F\hq \W$ is irreducible as well (see \cite[Ch.~II.4.3.A]{Kra84}).
\end{rem}

\begin{cor}\label{c:lunarichardson}
Consider the same notation as in the previous proposition. Let $y\in X_{\mathrm{pr}}$ and $b\in O(y)\cap X^H$. Consider the trivial linearisation on $V$ and $F$. Then $b\in F^s$.
\end{cor}

\proof
We start by proving that $a\in F^{s}$. First of all we have to show that $\W_a$ is the trivial group. This is true by the definition of $\W$. Secondly, we want to show that $\W(a)$ is closed in $F$. Assume that $\W(a)$ is not closed in $F$ and that $v\in \overline{\W(a)}\bsl \W(a)$. Since $G(a)=O(x)$ is closed in $V$, there exists a $g\in G$ such that $g.v=a$. If $g\in N_G(H)$, then there exists a $g'\in \W$ such that $g'.v=a$, which contradicts the assumption. Thus $g\notin N_G(H)$. However, in this case we have
$$a=g.v=g.(h.v)\neq h.(g.v)=h.a=a \text{, for }h\in H.$$
This is a contradiction, which shows that $v\in F^s$.

By \citep[Lem.~3.5]{LR79} we can deduce that $ O(y)\cap X^H\neq \emptyset$ and that $G_b$ is conjugated to $H$. Therefore, it is enough to show that $\W(b)$ is closed in $F$. Since $G(b)$ is closed in $V$, this follows from the same argumentation as above.
\eproof

\begin{ex}
Let $V=\CC^2$ and consider the following action of $G=\CC^*$ on $V$:
$$t.(z_1,\,z_2)=(z_1,\,t\cdot z_2),\text{ for }t\in \CC^*\text{ and }z_1,\,z_2\in V.$$
Using the trivial linearisation on $V$, every point is semi-stable (i.e. $V^{ss}=V$). Unfortunately, no point is stable, because the only closed orbits are $G(v)$ for points contained in the $(1,0)$-axis, $v\in \CC\cdot (1,0)$. These points however are all fixed by the group $G$ and thus $G_v=G$ is not finite, for all $v\in \CC\cdot (1,0)$. 

The quotient $X:=V\hq G$ is isomorphic to the space of complex numbers $\CC$. For every two points $x,\,x'\in X$, the unique closed orbits $O(x)=\{v\}$ and $O(x')=\{v'\}$ each consist of a single point contained in the $(1,0)$-axis and the stabiliser $G_v=G=G_{v'}$. Thus, $X_{\mathrm{pr}}=X$. Let $a\in O(0)$. Then $H=G$, $\W$ is the trivial group, $F=V^G=\CC\cdot (1,0)$ and $F\hq \W\cong X$. While the variety $X$ stays the same, the representation via $\W$ and $F$ is superior since $F^s=F$ and $F/\W$ is a geometric quotient. 
\end{ex}

\section{Some easy Extension Theorems and useful reductions}\label{ch.easyExtThm}

In this chapter we will present useful statements that will help us to simplify the proofs of the Extension Theorems in this paper. We then discuss special situations in which an Extension Theorem is easy to prove. For the convenience of the reader, let us recall the Extension Problem in the case of a normal variety:\\

\noindent\textbf{Extension Problem.}
Let $X$ be an $n$-dimensional normal variety and $\eta:\X\rightarrow X$ any resolution of $X$. For all $0\leq p\leq n$, we want to show that given a reflexive $p$-form $\s\in H^0\big(X,\,\Om_{X}^{[p]}\big)$ on $X$ the pull-back 
$$
\eta^*(\s)\in H^0\big(\X,\,\Om_{\X}^{p}\big)
$$
is a regular $p$-form on $\X$.\\

Recall that the assertion of the Extension Problem is local. This means that it is enough to show that $\eta^*(\s)$ has no pole near a fiber $\eta^{-1}(x)$, for all $x\in X$. The following result shows that for this purpose it is enough to consider an analytic neighbourhood of $x\in X$.

\begin{lem} \label{f:an}
In the same setting as above, let $U \subset X$ be an analytic open subset of $X$ with $x\in U$ and $\eta^{\mathrm{an}}$ the analytification of $\eta$. Assume that the form 
$$(\eta^*(\s))^{\mathrm{an}}_{\mid (\eta^{\mathrm{an}})^{-1}(U)}~~,$$
which is the restriction of the analytification of the pull-back of $\sigma$ to $(\eta^{\mathrm{an}})^{-1}(U)$, is a holomorphic form. Then $\eta^*(\s)$ has no pole near a fiber $\eta^{-1}(x)$.
\end{lem}

\proof
We want to show that a rational algebraic differential form on a smooth variety is regular if its analytification is a holomorphic form on the corresponding complex manifold. Using local coordinates, this follows from the fact that a rational algebraic function on a smooth variety is regular if its analytification is a holomorphic function on the corresponding complex manifold (see \cite[Book~3, p.~177]{Sha77b}).
\eproof

To avoid complicated notation when switching between the analytic and algebraic category we will not introduce the analytification of the relevant objects. However, we will use Lemma \ref{f:an} to restrict differential forms to analytic neighbourhoods when analysing their poles. 
\subsection{Independence of choice of resolution}\label{c:ind.of choice}

\begin{lem}\label{l:surj ext}
Let $Z$ be a normal variety and $Y$ a smooth variety. Consider a surjective morphism $f:Z\rightarrow Y$. Let $\s$ be a rational differential form (i.e. a rational section of the sheaf $\Om^{*}_Y$ of Kähler differential forms) on $Y$. Then $\s$ is regular if and only if the pull-back $f^*(\s)$ has no poles on $Z$.
\end{lem}

\proof
See \citep[Lem.~2]{Kem77}. 
\eproof

\begin{cor}[Independence of the choice of resolution] \label{c:choice of res}
Let $X$ be an $n$-dimensional normal variety and $\eta_1:\X_1\rightarrow X$ any resolution of $X$. Assume that for any reflexive $p$-form $\s\in H^0\big(X,\,\Om_{X}^{[p]}\big)$ ($0\leq p\leq n$) the pull-back
$$
\eta_1^*(\s)\in H^0\big(\X_1,\,\Om_{\X_1}^{p}\big)
$$
is a regular $p$-form on $\X_1$. Let $\eta_2:\X_2\rightarrow X$ be another resolution of $X$. Then
$$
\eta_2^*(\s)\in H^0\big(\X_2,\,\Om_{\X_2}^{p}\big)
$$
is a regular $p$-form on $\X_2$, too.
\end{cor}

\proof
There exists a resolution $\eta:\X\rightarrow X$ that dominates $\X_1$ and $\X_2$. In other words we can consider the following commutative diagram:
$$
	\xymatrix{
      \X	\ar[r]^{p_2}\ar[d]_{p_1} \ar[dr]^{\eta}	         &		 \X_2 \ar[d]^{\eta_2} \\
      \X_1	\ar[r]_{\eta_1}					 & 		~X .}
$$
Applying Lemma \ref{l:surj ext} to the maps $p_1$ and $p_2$ the result follows. 
\eproof

\subsection{Reductions for Extension Theorems for good quotients}\label{ch.red. ext}
In this chapter we want to analyse Extension Theorems for an arbitrary good quotient $X:=V\hq G$ of a smooth variety $V$ by a reductive group $G$ and present some useful reductions.

\begin{lem}[Reduction to affine GIT-quotients]\label{l:redtoaff}
To prove an Extension Theorem for a good quotient $X:=V\hq G$ (of a smooth variety $V$ by a reductive group $G$) it is enough to prove the Extension Theorem in the case where $V$ is an affine $G$-variety with affine GIT-quotient $X:=V\hq G$.
\end{lem}

\proof
Since proving the Extension Theorem is a local problem the result follows from the fact, that every good quotient can locally be described as an affine GIT-quotient (see Remark to Definition \ref{d:affineGIT}).
\eproof

\begin{lem}\label{l:reflexiveflat}
Let $g:Y\rightarrow X$ be a surjective flat morphism between normal varieties and $\mathcal F$ a coherent sheaf on $X$. Then the following assertions are true:

\enumeratel

\begin{enumerate}
	\item For the pull-back of the dual $\mathcal F^{\vee}$ of the sheaf $\mathcal F$ we have the isomorphism $g^*(\mathcal F^{\vee})\cong (g^*\mathcal F)^{\vee}$.
	\item The pull-back $g^*\mathcal F$ is reflexive if and only if $\mathcal F$ is reflexive.
\end{enumerate}
\end{lem}

\proof 
A proof of the first assertion can be found in \cite[1.8]{Har80}, where he explains that the problem is local and thus proves the following result for modules on rings instead:

Let $A$ be a noetherian ring, $M$, $N$ be $A$-modules, with $M$ finitely generated, and let $A\rightarrow B$ be a flat surjective ring homomorphism (for the first assertion we do not need $g$ to be surjective). Then the natural map
$$
\Hom_A(M,N)\otimes_A B\rightarrow \Hom_B\big(M\otimes_A B, N\otimes_A B\big)
$$
is an isomorphism (see \cite[1.8]{Har80}).

To prove the second assertion we have to show that the natural map $\mathcal F\rightarrow \mathcal F^{\vee\vee}$ is an isomorphism if and only if the natural map $g^*\mathcal F\rightarrow (g^*\mathcal F)^{\vee\vee}$ is an isomorphism. Using \ref{l:reflexiveflat}.$1$ we have $(g^*\mathcal F)^{\vee \vee}\cong g^*(\mathcal F^{\vee \vee})$. Since proving the equivalence is a local problem, too, it can be deduced from the following claim:

\noindent\textbf{Claim 1.}
Consider the same notation as above. Then the natural map 
$$ M\rightarrow \Hom_A(\Hom_A(M,N),N)$$ 
is an isomorphism if and only if the map 
$$M\otimes_A B\rightarrow\Hom_A(\Hom_A(M,N),N)\otimes_A B$$ 
is an isomorphism.

\proofof{Claim 1}
Since $g$ is flat and surjective, $B$ is a faithfully flat $A$-module. Thus tensoring a sequence of $A$-modules with $B$ over $A$ produces an exact sequence if and only if the original sequence was exact. 
\eproof

\begin{lem}\label{l:etaleEXT}
Let $g:Y\rightarrow X$ be a surjective \'{e}tale morphism between normal varieties. Then the Extension Theorem for $X$ is true if and only if it is true for $Y$. 
\end{lem}

\proof 
Let $\eta_X:\X\rightarrow X$ be a resolution of singularities of $X$. Consider the following commutative diagram
$$
\xymatrix{
      \widetilde Y	\ar[r]^{\widetilde g}\ar[d]_{{\eta}_Y}          &		 \X \ar[d]^{{\eta}_X} \\
       Y	\ar[r]_{g}					 & 		X,}
$$
where $\widetilde Y$ is the fiber product $\widetilde Y=Y\times_X \X$. Thus, the morphism $\widetilde g$ is \'{e}tale (see \cite[Prop.~4.3.(iii)]{Dre12}) and $\widetilde Y$ is smooth (see \cite[Prop.~4.3.(vii)]{Dre12}). Since $\eta_X$ is a resolution, the morphism $\eta_Y:\widetilde Y\rightarrow Y$ is a resolution of singularities of $Y$. Let $0\leq p\leq n=\dim(X)=\dim(Y)$. Then we want to show that $(\eta_X)_*\Om^p_{\X}$ is reflexive if and only if the sheaf $(\eta_Y)_*\Om^p_{\widetilde Y}$ is reflexive (this is enough since the Extension Theorem is independent of the choice of resolution (see Corollary \ref{c:choice of res})). 

Since for $p=0$ both sheaves are reflexive (see Remark to Proposition \ref{nforms}), we may assume that $p>0$. Since $\widetilde g$ is \'{e}tale the sheaf of relative differential forms $\Om^p_{\widetilde{Y}/\X}=0$ (see \cite[Prop.~4.2]{Dre12}) and thus $\Om^p_{\widetilde Y}\cong (\widetilde{g})^*\Om^p_{\X}$. Since $g$ is \'{e}tale it is flat in particular. Using the cohomology-and-base-change (see \cite[III:~Prop.~9.3]{Har77}) we get 
$$
(\eta_Y)_*\Om^p_{\widetilde Y}\cong g^*\big((\eta_X)_*\Om^p_{\X}\big)
$$
and the result follows from Lemma \ref{l:reflexiveflat}.
\eproof

\begin{lem}\label{l:redtoV}
Let $G$ be a reductive group and $V$ a smooth affine $G$-variety with affine GIT-quotient $X:=V\hq G$ and quotient map $\pi:V\rightarrow X$. Let $x\in X$ be any point in $X$. Then there exists a vector space $W$ and a reductive subgroup $H\subset G$ acting linearly on $W$, such that the Extension Theorem for the quotient $W\hq H$ implies the Extension Theorem for $X$ in a neighbourhood of $x\in X$.
\end{lem}

\proof
Let $x\in X$ and $v\in\pi^{-1}(x)$ be a point with closed $G$-orbit $G(v)$ and therefore reductive stabiliser $G_v$ (see Lemma \ref{l:Mats}). By Theorem \ref{t:LST} there exists a locally closed affine $G_v$-invariant sub-variety $S\subset V$ with $v\in S$ and an open saturated subset $U\subset V$ such that
$$
f:S\hq G_v\rightarrow U\hq G
$$
is a surjective \'{e}tale morphism.
Since $v\in V$ is smooth we also get a linear action of $G_v$ on the tangent space $T_vS$ to the slice $S$ at $v\in S$ and an open saturated subset $U'\subset T_vS$, such that 
$$
g:S\hq G_v\rightarrow U'\hq G_v
$$
is a surjective \'{e}tale morphism. Set $W:=T_vS$ and $H:=G_v$. Then the result follows by applying Lemma \ref{l:etaleEXT} to $f$ and $g$. 
\eproof

\begin{cor}[Reduction to quotients of vector spaces]\label{p:redtoV}
To prove an Extension Theorem for an arbitrary good quotient $X:=V\hq G$ of a smooth variety by a reductive group it is enough to prove the Extension Theorem in the case where $V$ a vector space with linear $G$-action and affine GIT-quotient $X:=V\hq G$.
\end{cor}

\proof
Since proving the Extension Theorem is a local problem we can use Lemma \ref{l:redtoaff} to reduce to the case  where $V$ is an affine $G$-variety with affine GIT-quotient $X:=V\hq G$. Using the fact that proving the Extension Theorem for the affine GIT-quotient is a local problem as well, the result follows from Lemma \ref{l:redtoV}.
\eproof

\subsection{Finite quotient singularities}

A special type of singularities that can arise on good quotients of smooth varieties are finite quotient singularities. In this chapter we will discuss these singularities and prove an Extension Theorem for good quotients with finite quotient singularities. 

\begin{defn}[Finite quotient singularity]
Let $X$ be a normal $n$-dimensional variety. Then a point $x\in X$ is called \emph{finite quotient singularity (f.q.s)} if there exists an analytic neighbourhood $U$ of $x$ and a finite group $\Gamma$ acting linearly on some $\CC^n$, such that $U$ is biholomorphic to an open neighbourhood $U'\subset\CC^n/\Gamma$ of the image of $0\in \CC^n$ in the quotient $\CC^n/\Gamma$.
\end{defn}

\begin{ex}
The following two examples will be revisited in Chapter \ref{subs:cutting down}.
\enumeraten
\begin{enumerate}
\item If $X$ has klt singularities, then there exists a closed subset $Z\subset X$ with $\codim_X(Z)\geq 3$, such that $X\bsl Z$ has only finite quotient singularities (see \citep[Prop.~9.4]{GKKP11}).
\item Let $G$ be a reductive group acting linearly on $\CC^n$ such that the quotient $X:=\CC^n\hq G$ is a surface with an isolated singularity $t\in X$. Then $t\in X$ is a finite quotient singularity (see \citep[Cor.~2]{Gur}).
\end{enumerate}
\end{ex}

We will discuss the following example separately, because it plays an important role in Chapter \ref{ch. Kirwan}.  

\begin{ex} \label{e:fqs stable}
Let $V$ be an affine or projective smooth $G$-variety admitting a GIT-quotient $X:=V^{ss}\hq G$ (with respect to some $G$-linearisation). Let $\pi:V^{ss}\rightarrow X$ be the induced quotient map and assume that $V^s=V^{ss}$. Then $\pi:V^{ss}\rightarrow X$ is a geometric quotient and $X$ has finite quotient singularities due to the following argument:

Let $x\in X$ be an arbitrary point. By shrinking $X$ to a neighbourhood of $x\in X$, we may assume that $V$ is affine and $X$ is an affine GIT-quotient (see Remark to Definition \ref{d:affineGIT}). Let $v\in \pi^{-1}(x)$ be a point with closed orbit $G(v)$ and therefore reductive stabiliser $H:=G_v$ (see Lemma \ref{l:Mats}). Then (by Corollary \ref{l:gurjar}) $H$ acts linearly on some $\CC^n$ such that the analytic local ring of $X$ in $x$ is isomorphic to the analytic local ring of $\CC^n\hq H$ at the image of $0\in \CC^n$ in $\CC^n\hq H$. Since $v\in V^s$, we know that $H$ is finite. Thus $x\in X$ is a finite quotient singularity.
\end{ex}

\begin{prop}\label{p:ExtThm fqs}
Let $G$ be a reductive group and $V$ a smooth $G$-variety admitting a good quotient $X:=V\hq G$. Assume that $X$ has f.q.s. Let $\eta:\X\rightarrow X$ be a resolution and let $\s\in H^0\big(X,\,\Om_{X}^{[p]}\big)$ be a reflexive $p$-form on $X$, for $0\leq p\leq \dim(X)$. Then 
$$
\eta^*(\s)\in H^0\big(\X,\,\Om_{\X}^{p}\big).
$$
\end{prop}

\proof
The result can be deduced from the fact that f.q.s. are klt (see \citep[Prop.~5.20]{KM98}) by using \citep[Thm.~1.4]{GKKP11}.
\eproof
 
\subsection{Good quotients that are rational Gorenstein}

\begin{prop}\label{p:ExtThm Gorenstein}
Let $G$ be a reductive group and $V$ a smooth $G$-variety admitting a good quotient $X:=V\hq G$. Assume that the canonical divisor $K_X$ is Cartier. Let $\eta:\X\rightarrow X$ be a resolution and $\s\in H^0\big(X,\,\Om_{X}^{[p]}\big)$ a reflexive $p$-form on $X$, for $0\leq p\leq \dim(X)$. Then 
$$
\eta^*(\s)\in H^0\big(\X,\,\Om_{\X}^{p}\big).
$$
\end{prop}

\proof
Since $X$ has rational singularities, due to \citep[Cor.~5.24]{KM98} $K_X$ Cartier implies that $X$ has canonical singularities (see \citep[Def.~2.11]{KM98}) and therefore klt singularities. Thus the assertion follows from \citep[Thm.~4.1]{GKKP11}.
\eproof

It is important to note, that not every GIT-quotient $X=V\hq G$ (of a smooth variety $V$ by a reductive group $G$) is Gorenstein. An example can be constructed using the quotient construction of a toric variety.

\begin{ex}
Let $X$ be an affine toric variety (see \cite[Def.~1.1.3]{CLS11}). Then Cox-Little-Schenck describe a quotient construction of $X$ in \citep[§~5.1]{CLS11}. As a consequence, we can think of $X$ as a quotient $V\hq G$ where $G$ is a reductive group and $V$ is a smooth affine $G$-variety. In order to find a GIT-quotient that is not Gorenstein we thus only have to give an example of a toric variety that is not Gorenstein.

In \citep[§1.2]{CLS11} one can find the definition of a cone $c$ and a discussion on how it defines an affine toric variety $U_{c}$. In the examples \cite[Ex.~1.2.22, Ex.~4.1.4, Ex.~8.2.13]{CLS11} Cox-Little-Schenck discuss the affine toric variety $U_{c}$ that is induced by the cone $c:=\cone(de_1-e_2,e_2)\subset \RR^2$, where $e_1,e_2$ are the standard basis of $\RR^2$ and $d\in \NN$ is a positive integer. Using the divisors $D_1$, $D_2$ corresponding to the rays of $c$ (see Orbit-Cone Correspondence in \cite[Thm.~3.2.6]{CLS11}) one can show that $K_{U_{c}}=-D_1-D_2$ (\cite[Thm.~8.2.3]{CLS11}). Using \citep[Prop.~8.2.12]{CLS11} Cox-Litte-Schenk show that $U_{c}$ is Gorenstein if and only if $d\leq 2$. 
\end{ex}

\section{Preparation for Theorem \ref{t:2}}

\subsection{Hodge-theoretic method}\label{Ch.MHS}
In this chapter we will present the Hodge-theoretic approach to the Extension Theorem. Most of the ideas are based on \citep{vSS85} and \citep{Nam01}. Following the guidelines of these papers we will work in the following setting:

\vspace{10pt}
\fboxsep3mm
\noindent\fbox{
\begin{minipage}[c]{0.945\textwidth}
\begin{set}
Let $X$ be a Stein open subset of an algebraic variety with rational singularities and $x\in X$ a point in $X$. Let $\eta:\X\rightarrow X$ be a resolution of $X$ such that $E:=(\eta^{-1}(x))_{\mathrm{red}}$ is an snc divisor on $\X$ with support $\eta^{-1}(x)$. 
\end{set}
\end{minipage}
}
\vspace{10pt}

\begin{rem}[Working in the analytic category]
For the first part of this chapter we do not have to assume that $X$ is an analytic space. In Lemma \ref{MHSdelta} however it is necessary to consider a small analytic neighbourhood of $x\in X$ and in the proof of Corollary \ref{surjection} we work in this setting, too. 
\end{rem}

The main goal of this chapter is to show that the injection
$$
\iota: H^0\big(\X,\,\Om_{\X}^{p} \big) \rightarrow H^0\big(\X,\,\Om_{\X}^{p}\big( \log E \big)\big)
$$
is an isomorphism for all values $0\leq p\leq \dim(X)$. To do so we consider the exact sequence (see \citep[p.~8]{Nam01})
$$
0\rightarrow \Om_{\X}^{p}/ \Om_{\X}^{p}\big( \log E \big)\big( - E \big)\rightarrow \Om_{\X}^{p}\big( \log E \big)/ \Om_{\X}^{p}\big( \log E \big)\big( - E \big) \rightarrow \Om_{\X}^{p}\big( \log E \big)/ \Om_{\X}^{p}\rightarrow 0
$$
of sheaves on $E$. We will refer to this sequence by (+)\label{sequence}. 

\begin{lem}\label{deltainjective}
Let $\delta: H^0\big(\X,\,\Om_{\X}^{p}\big( \log E \big)/ \Om_{\X}^{p} \big) \rightarrow H^1\big(\X,\,\Om_{\X}^{p}/\Om_{\X}^{p} \big( \log E \big)\big( - E \big) \big)$ be the boundary map in the long exact cohomology sequence induced by (+). If $\delta$ is injective for some value of $p\leq \dim(\X)$, then $\iota$ is an isomorphism for the same value of $p$. 
\end{lem}

\proof
The following argument goes back to \citep[Thm.~1.3]{vSS85}. Consider the long exact sequence
$$
0\rightarrow  H^0\big(\X,\,\Om_{\X}^{p} \big) \xrightarrow{\iota}  H^0\big(\X,\,\Om_{\X}^{p}\big( \log E \big)\big)\rightarrow H^0\big(E,\,\Om_{\X}^{p}\big( \log E \big)/\Om_{\X}^{p}\big) \xrightarrow{\gamma}  H^1\big(\X,\,\Om_{\X}^{p} \big)\rightarrow \dots .
$$
Then by \citep[p.~99]{vSS85} $\delta$ can be understood as the composition 
$$
H^0\big(E,\,\Om_{\X}^{p}\big( \log E \big)/\Om_{\X}^{p}\big) \xrightarrow{\gamma}  H^1\big(\X,\,\Om_{\X}^{p} \big)\rightarrow H^1\big(E,\,\Om_{\X}^{p}/\Om_{\X}^{p} \big( \log E \big)\big( - E \big) \big).
$$
Thus, we get the following diagram:

$$
\xymatrix{
0\ar[d] & \\
H^0\big(\X,\,\Om_{\X}^{p} \big) \ar[d]^{\iota}  & \\
H^0\big(\X,\,\Om_{\X}^{p}\big( \log E \big)\big)\ar[d] & \\
H^0\big(E,\,\Om_{\X}^{p}\big( \log E \big)/\Om_{\X}^{p}\big) \ar[d]^{\gamma} \ar[dr]^{\delta}& \\
H^1\big(\X,\,\Om_{\X}^{p} \big)\ar[d] \ar[r] & H^1\big(E,\,\Om_{\X}^{p}/\Om_{\X}^{p} \big( \log E \big)\big( - E \big) \big) \\
\dots & .
}
$$
Since $\delta$ is injective the map $\gamma$ is injective as well. Since the vertical sequence is exact, $\iota$ must be surjective. 
\eproof

\begin{rem}
The horizontal map in the last diagram is induced by the exact sequence
$$
0\rightarrow \Om_{\X}^{p} \big( \log E \big)\big( - E \big)\rightarrow \Om_{\X}^{p}\rightarrow  \hat{\Om}_{E}^{p} \rightarrow 0.
$$ 
of sheaves on $E$. The sheaf $\hat{\Om}_{E}^{p}:=\Om_{\X}^{p}/\Om_{\X}^{p} \big( \log E \big)\big( - E \big)$ defined by this sequence is called \emph{sheaf of torsion free $p$-forms on $E$}. A prove of its properties can be found in \citep[Part~I]{Keb13}.
\end{rem}

\begin{rem}
We would like to recall the fact that after shrinking $X$ to a small analytic neighbourhood of $x\in X$, we get an isomorphism of the cohomology groups $H^k(\X,\,\CC)\cong H^k(E,\,\CC)$ (see \citep[p.~7]{Nam01}). We will now work in this setting.
\end{rem}

\begin{lem}\label{MHSdelta}
The cohomology groups $H^k(\X,\,\CC)$ and $H_{E}^k(\X,\,\CC)$ carry mixed Hodge structures (MHS) with filtrations $F$ and $W$, such that the boundary morphism 

$$ 
\delta: H^0\big(E,\,\Om_{\X}^{p}\big( \log E \big)/ \Om_{\X}^{p} \big) \rightarrow H^1\big(E,\,\Om_{\X}^{p}/\Om_{\X}^{p} \big( \log E \big)\big( - E \big) \big)
$$
can be interpreted as the map 
$$
Gr_F^p H_{E}^{p+1}(\X,\,\CC)\rightarrow Gr_F^p H^{p+1}(\X,\,\CC)
$$
for all $0\leq p\leq \dim(\X)$.
\end{lem}

\proof
The proof of this lemma follows from \citep[pp.~7-9]{Nam01}.
\eproof

\begin{rem}
A good description of the MHS on $H^k(\X,\,\CC)$ and $H_{E}^k(\X,\,\CC)$ can be found in \citep[(1.5), (1.6)]{Ste83}. 
\end{rem}

By combining the previous two lemmas the problem of $\iota$ being an isomorphism reduces to the following corollary. 

\begin{cor} \label{surjection}
Let $U:=\X\bsl E$ be the complement of $E$ in $\X$. If the morphism of cohomology groups $ \alpha: H^{p}(\X,\,\CC)\rightarrow H^{p}(U,\,\CC) $ is a surjection for some value of $p\leq \dim(\X)$, then $\iota$ is an isomorphism for the same value of $p$. 
\end{cor}

\proof
Consider the exact local cohomology sequence (see \citep[p.~7]{Nam01}) 
$$
\dots \rightarrow  H^{p}(\X,\,\CC) \xrightarrow{\alpha} H^{p}(U,\,\CC)\rightarrow H_{E}^{p+1}(\X,\,\CC)\xrightarrow{\beta} H^{p+1}(\X,\,\CC) \rightarrow \dots .
$$
Since $\alpha$ is surjective the exactness implies that $\beta$ is injective. The map $\beta$ however is a morphism of mixed Hodge structures. Thus, $Gr_F^p H_{E}^{p+1}(\X,\,\CC)\rightarrow Gr_F^p H^{p+1}(\X,\,\CC)$ is injecive. The statement of the corollary now follows by using Lemma \ref{MHSdelta} and Lemma \ref{deltainjective}.
\eproof

\subsection{Two Extension Theorems} 
In this chapter we will recall two important Extension Theorems that were proven using the method presented in the previous chapter. In both cases the result follows from Corollary \ref{surjection}. 

The first result on spaces with isolated singularities was proven by van Straten-Steenbrink in \citep[Thm.~1.3]{vSS85}. In their paper they work on an arbitrary contractible Stein space with an isolated singularity. Thus, their proof, although still using Corollary \ref{surjection}, slightly differs from our proof. The case $p\leq n$ can be found in \citep[Thm.~1.3]{vSS85}, whereas the case $p=n$ follows from \citep[Cor.~1.4]{vSS85}.

\begin{cor}[Extension Theorem for isolated rational singularities]
Let $X$ be an affine variety with rational singularities, $\dim(X):=n\geq 2$ and $x\in X$ the only singular point of $X$. Let $\eta:\X\rightarrow X$ be a log-resolution of $X$ such that $E:=\exc(\eta)$ is an snc divisor on $\X$. For $0\leq p \leq n$, let $\s\in H^0\big(X,\,\Om_{X}^{[p]}\big)$ be a reflexive $p$-form on $X$. Then 
$$
\eta^*(\s)\in H^0\big(\X,\,\Om_{\X}^{p}\big).
$$

\end{cor}

\proof
We have to consider two cases. If $p=n$ the result follows from the Remark to Proposition \ref{nforms}. Assume that $p< n$. By Proposition \ref{logpoles} the pull-back $\eta^*(\s)\in H^0\big(\X,\,\Om_{\X}^{p}\big( \log E \big)\big)$ is a regular $p$-form on $\X$ with logarithmic poles along $E$. Let $E'$ be the divisor that contains all components of $E$ that are mapped to $x\in X$ via $\eta$. Since $E$ is a reduced snc divisor, $E'$ is a reduced snc divisor as well. Note that $X$ is smooth outside of $x\in X$. Thus $\eta^*(\s)\in H^0\big(\X,\,\Om_{\X}^{p}\big( \log E' \big)\big)$. We only have to show, that $\eta^*(\s)$ has no logarithmic poles along $E'$. Since proving the Extension Theorem is a local problem by Lemma \ref{f:an} we can replace $X$ by a small Stein open subset of an algebraic variety such that we get an isomorphism $H^k(\X,\,\CC)\cong H^k(E',\,\CC)$ (see Remark before Lemma \ref{MHSdelta}). We want to use Corollary \ref{surjection} and thus have to show that the restriction map
$$
\alpha: H^{p}(\X,\,\CC)\rightarrow H^{p}(U,\,\CC) 
$$
is a surjection for $p< n$, where $U:=\X\bsl E'$ is the complement of $E'$ in $\X$. This statement is a corollary of the decomposition theorem in intersection cohomology and goes back to Goreski and MacPherson (see \cite[Thm.~1.11]{Ste83}). 
\eproof

The second result is due to Namikawa \citep[Lem.~2]{Nam01}, who works in the setting of Chapter \ref{Ch.MHS}.

\begin{cor}\label{c:12forms}
Let $X$ be a Stein open subset of an algebraic variety with rational singularities of dimension $n\geq 3$ and $x\in X$ a point in $X$. Let $\eta:\X\rightarrow X$ be a resolution of $X$ such that $E:=(\eta^{-1}(x))_{\mathrm{red}}$ is an snc divisor on $\X$. Then 
$$
H^0\big(\X,\,\Om_{\X}^{p}\big)\rightarrow H^0\big(\X,\,\Om_{\X}^{p}\big( \log E \big)\big)
$$
is an isomorphism for $p\in\{1,2\}$.
\end{cor}

\proof

We once again we want to use Corollary \ref{surjection} and thus have to show that 
$$
\alpha: H^{p}(\X,\,\CC)\rightarrow H^{p}(U,\,\CC) 
$$
is a surjection for $p \in\{1,2\}$, where $U:=\X\bsl E$ is the complement of $E$ in $\X$. Namikawa proves this in \citep[Lem.~2]{Nam01} by exploiting the geometric structure of $H^{p-1}(\X,\,\mathcal O_{X}^*)$.
\eproof

The following remark will be useful for the proof of Theorem \ref{t:2} in Chapter \ref{ch. proofof1.2}.

\begin{rem}
Consider the exact sequence (+)
$$
0\rightarrow \Om_{\X}^{p}/ \Om_{\X}^{p}\big( \log E \big)\big( - E \big)\rightarrow \Om_{\X}^{p}\big( \log E \big)/ \Om_{\X}^{p}\big( \log E \big)\big( - E \big) \rightarrow \Om_{\X}^{p}\big( \log E \big)/ \Om_{\X}^{p}\rightarrow 0
$$
of sheaves on $E$ we saw in Chapter \ref{Ch.MHS}. Because $\alpha$ is surjective for $p \in\{1,2\}$ by Lemma \ref{MHSdelta} and Lemma \ref{deltainjective} we can deduce that the boundary map
$$
\delta: H^0\big(E,\,\Om_{\X}^{p}\big( \log E \big)/ \Om_{\X}^{p} \big) \rightarrow H^1\big(E,\,\Om_{\X}^{p}/\Om_{\X}^{p} \big( \log E \big)\big( - E \big) \big)
$$
in the long exact sequence is injective for these values of $p$. Thus
$$
H^0\big(E,\,\Om_{\X}^{p}/ \Om_{\X}^{p}\big( \log E \big)\big( - E \big)\big) \rightarrow H^0\big(E,\,\Om_{\X}^{p}\big( \log E \big)/\Om_{\X}^{p} \big( \log E \big)\big( - E \big) \big)
$$
is surjective for $p \in\{1,2\}$.
\end{rem}

\begin{cor}[Extension Theorem for $1$-, $2$-forms on $(X,x)$]\label{c:12formsaffine}
Let $X$ be an affine variety with rational singularities of dimension $n\geq 3$ and $x\in X$ a point in $X$. Let $\eta:\X\rightarrow X$ be a resolution of $X$ such that $E:=(\eta^{-1}(x))_{\mathrm{red}}$ is an snc divisor on $\X$. For $0\leq p\leq n$, let $\s\in H^0\big(X,\,\Om_{X}^{[p]}\big)$ be a reflexive $p$-form on $X$ such that the pull-back $\eta^*(\s)\in H^0\big(\X,\,\Om_{\X}^{p}\big( \log E \big)\big)$ is a regular $p$-form on $\X$ with a potential logarithmic poles along $E$. Then $\eta^*(\s)$ has no pole along $E$, for all values $0\leq p \leq n$.
\end{cor}

\proof
Since proving the Extension Theorem is a local problem we can restrict ourselves to a Stein open subset of $X$ such that the conditions of Chapter \ref{Ch.MHS} are fulfilled. By Lemma \ref{f:an} the result now follows directly from the previous corollary.
\eproof

\subsection{Cutting down}\label{subs:cutting down}
In this chapter we will show that affine GIT-quotients are stable under general hyperplane sections. We later use the results of this chapter to show that GIT-quotient have finite quotient singularities in codimension $2$ and to prove Theorem \ref{t:2}. 

\begin{lem}\label{l:cutting down}
Let $G$ be a reductive group and $V$ a smooth affine $G$-variety. Assume that the affine GIT-quotient $X:=V\hq G$ has dimension $dim(X)\geq 2$ and let $\pi:V\rightarrow X$ be the quotient map. Let $H\in |\mathscr L |$ be a general element of an ample basepoint-free linear system corresponding to $\mathscr L\in Pic(X)$. Then the following statements hold.
\enumeratel
\begin{enumerate}
 \item The divisor $H$ is irreducible and normal.
 \item If $H$ is smooth, then $X$ is smooth along $H$.
 \item If $\eta:\X\rightarrow X$ is a (strong) log-resolution and $\widetilde H:=\pi^{-1}(H)$, then the restriction $\eta_{\mid H}:\widetilde H\rightarrow H$ is a (strong) log-resolution with exceptional set $\exc(\eta_{\mid H})=\exc(\eta)\cap H$.
 \item The preimage $\pi^{-1}(H)$ is a smooth affine G-invariant hyperplane in $V$ and $H$ is the affine GIT-quotient $H=\pi^{-1}(H)\hq G$.
\end{enumerate}
\end{lem}

\proof
The first three statements can be found in \citep[Lem.~2.22, Lem.~2.23]{GKKP11}. Assertion \ref{l:cutting down}.4 can be deduced from the following observations. If $\pi(v)=x$ for some $x\in H$ and $v\in V$ then $\pi(g.v)=x$ for all $g\in G$. Assume that for $1\leq l\leq n$ $f_1,\,\dots\, f_l\in \CC[X]$ are the defining polynomials of $H=\{x\in X\mid f_1(x)=\dots=f_l(x)=0\}$. Then $\pi^{-1}(H)=\{v\in V\mid f_1\circ\pi(v)=\dots=f_l\circ\pi(v)=0 \}$. Thus $\pi^{-1}(H)$ is a $G$-invariant hyperplane on V. This hyperplane is smooth by Bertini's Theorem (c.f. \citep[II.8.18]{Har77}) and affine since it is the fiber product in the pull-back diagram
$$
\xymatrix{
\pi^{-1}(H)\ar[r]\ar[d] & V\ar[d] \\
H\ar[r] & X,
}
$$
where $H$, $X$ and $V$ are affine. 
\eproof

\part{The 3-dimensional case}
This part is dedicated to the proof of Theorem \ref{t:2}. Recall that it is enough to prove Theorem \ref{t:2} in the case where the good quotient $X$ is an affine GIT-quotient. We will first show that GIT-quotients have finite quotient singularities in codimension $2$. Then we prove Theorem \ref{t:2} by adjusting the proof of \citep[Prop.~3]{Nam01} to $1$- and $2$-forms on affine GIT-quotients.

\section{Extension in codimension 2}

\subsection{Singularities of GIT-quotients in codimension $2$}\label{ch:singcodim2}
In this chapter we will generalise a result by Gurjar. He showed that $2$-dimensional affine GIT-quotients have finite quotient singularities. We will use a standard cutting down technique presented by Greb-Kebekus-Kovàcs-Peternell in the case of klt-pairs (see \citep[Ch.~9.C]{GKKP11}) to show that GIT-quotients have finite quotient singularities in codimension two. 

\begin{prop}[c.f. {\citep[Cor.~2]{Gur}}] \label{c:dim2}
Let $G$ be a reductive group and $V$ a smooth affine $G$-variety with affine GIT-quotient $X:=V\hq G$. Assume that $dim(X)=2$. Then $X$ has at most finite quotient singularities.
\end{prop}

\proof
By Lemma \ref{l:gurjar} we can assume that $X=\CC^n\hq H$, where $H$ is a reductive group acting linearly on $\CC^n$. Thus, Lemma \ref{l:goodC} asserts that $X$ has a good $\CC^*$-action. Therefore, \citep{Pin77} has shown that $X$ is locally isomorphic to $\CC^2\hq \Gamma$ where $\Gamma$ is a finite group acting linearly on $\CC^2$. \eproof

\begin{prop}[Affine GIT-quotients have finite quotient singularities in codimension 2]\label{p:codim2}
Let $G$ be a reductive group and $V$ a smooth affine $G$-variety with affine GIT quotient $X:=V\hq G$ and quotient map $\pi:V\rightarrow X$. Then there exists a closed subset $Z\subset X$ with $\codim_X(Z)\geq 3$ such that $X\bsl Z$ has finite quotient singularities. 
\end{prop}

\begin{rem}
The proposition of course is true if $dim(X)=0$ or $dim(X)=1$, because in both cases $X$ is smooth. The case $dim(X)=2$ is covered in Proposition \ref{c:dim2} and was proven by Gurjar. Due to this, we will prove Proposition \ref{p:codim2} only for $dim(X)\geq 3$.
\end{rem}

\proof (c.f. \citep[Prop.~9.4]{GKKP11})
We basically follow the proof presented in \citep[Ch.~9.C]{GKKP11} and fill in additional steps where they are needed. We start by reducing the problem to a more simple case:\\

\noindent\textbf{Step 1:}
Recall that $X$ is a normal variety and the singular locus $T:=X_{\mathrm{sing}}$ has codimension at least $2$. We can find a closed subset $T'\subset T$ such that every irreducible component of $T\bsl T'$ has codimension $2$ and $\codim_X(T')\geq 3$. The assertion of Proposition \ref{p:codim2} is local on $X$. Thus, we may assume that $T\subset X$ is irreducible with $\codim_X(T)=2$. 

By \citep[Prop.~2.26]{GKKP11} there exists an open set $X^0\subset X$ such that $T^0:=T\cap X^0$ is non-empty and a diagram
$$
\xymatrix{
\mathcal{X}^0 \ar[rr]^{\gamma}_{finite, \ etale} \ar[d]_{\F} & & X^0\\
S^0                                                   & &
}
$$
such that the restriction of $\F$ to any connected component of $\gamma^{-1}(T^0)$ is an isomorphism. The subspace $X\bsl (X^0\cup T)\subset X$ is smooth and $\codim_X(T\bsl T^0)\geq 3$. Consequently, it is sufficient to prove Proposition \ref{p:codim2} for points contained in $X^0$. Moreover, since the assertion in the Proposition \ref{p:codim2} is local in the analytic topology, it suffices to prove it for the variety $\mathcal{X}^0$ instead, even after removing all but one component of $\gamma^{-1}(T^0)$. In conclusion we may assume the following:

\begin{center}
\fbox{
\begin{minipage}[c]{0.8\textwidth}
There exists a surjective morphism $\F:X\rightarrow S$ with connected fibers, such that the restriction $\F_{\mid T}:T\rightarrow S$ is an isomorphism.
\end{minipage}
}
\end{center}

\noindent\textbf{Step 2:}
We now want to reduce our problem to a case where the fibers of the morphism $\F$ are surfaces. Let $S^0\subset S$ be a Zariski-open, dense subset. Then $X$ is smooth at all points of $X\bsl (\F^{-1}(S^0)\cup T)$ and $\codim_X(T\bsl\F^{-1}(S^0))\geq 3$. As above, it is sufficient to prove Proposition \ref{p:codim2} for the open set $\F^{-1}(S^0)\subset X$ only. Consider the following diagram
$$
\xymatrix{
V \ar[r]^{\pi} \ar[dr]_{\lambda} & X \ar[d]^{\F} & \\
& S &.
}
$$
For each $s\in S$ we denote the fibers by $X_s:=\F^{-1}(s)$ and $V_s:=\pi^{-1}(X_s)=\lambda^{-1}(s)$. Then the Generic Flatness Lemma \citep[Lem.~5.12]{FGI05} and the Lemma \ref{l:cutting down} allow us to assume the following:

\begin{center}
\fbox{
\begin{minipage}[c]{0.8\textwidth}
The morphism $\F$ is flat. Given any point $s\in S$ the preimage $\Xs:=\F^{-1}(s)$ is a normal surface, and moreover an affine GIT-quotient $X_s=V_s\hq G$, where $V_s$ is a smooth variety. If $t_s\in T$ is the unique point that maps to $s\in S$, then $\Xs$ is smooth away from $t_s\in \Xs$. Using Proposition \ref{c:dim2}, $X_s$ has only finite quotient singularities. In particular, $X_s$ is klt [KM98; 5.20]. 
\end{minipage}
}
\end{center}

\noindent\textbf{Step 3:}
We are now in the situation of \citep[9.8]{GKKP11} and can adopt their proof from here on. Proposition \ref{p:codim2} follows from the argument in \citep[9.C.2 and 9.C.3]{GKKP11}.
\eproof

\begin{rem}
Since the assertion of the previous proposition is local on $X$ it holds for good quotients as well. 
\end{rem}

\subsection{Two Extension Theorems}

\begin{cor}[Extension Theorem for GIT-quotients in codimension $2$]\label{c:codim2}
Let $X$ and $Z$ be as in Proposition \ref{p:codim2}. For $0\leq p\leq dim(X)$ let $\s\in H^0\big(X\bsl Z\, , \, \Omega_{X}^{[p]}\big)$ be a reflexive $p$-form on $X\bsl Z$ and $\eta:\X\rightarrow X$ a log-resolution. Then $\eta^*\s\in H^0\big(\eta^{-1}(X\bsl Z)\, , \,\Omega_{\X}^{p}\big)$.
\end{cor}

\proof This follows from the Extension Theorem for varieties with finite quotient singularities, Proposition \ref{p:ExtThm fqs}. \eproof

\begin{rem}
Similarly to the remark to Proposition \ref{p:codim2} this Extension Theorem holds for good quotients as well. 
\end{rem}

At this point we can already prove Theorem \ref{t:1} for $dim(X)=3$ using the previous corollary together with Corollary \ref{c:12formsaffine}.

\begin{cor}[Extension Theorem for good quotients in dimension 3]\label{c:Extdim3}
Let $G$ be a reductive group and $V$ a smooth $G$-variety admitting a good quotient $X:=V\hq G$ of dimension $dim(X)=3$. Let $\s\in H^0\big(X\, , \, \Omega_X^{[p]}\big)$ be a reflexive $p$-form, $0\leq p\leq dim(X)$, and $\eta:\X\rightarrow X$ a log-resolution. Then 
$$
\eta^*\s\in H^0\big(\X\, , \,\Omega_{\X}^{p}\big).
$$

\end{cor}

\proof
Without loss of generality $V$ is affine and $X$ is an affine GIT-quotient. The proposition is clear for $p=0$ and $p=n$ by the remark to Proposition \ref{nforms}. Thus we only have to check the assertion for $1\leq p\leq 2$. 

Let $\s\in H^0\big(X\, , \, \Omega_X^{[p]}\big)$ be a reflexive $p$-form, for $1\leq p\leq 2$. By Proposition \ref{logpoles}, $\eta^*\s\in H^0\big(\X\, , \,\Omega_{\X}^{p}(\slog E)\big)$, where $E:=\exc(\eta)$ is the exceptional divisor. Thus, we only have to show that $\s$ extends on divisors contained in $E$. Let $E'\subset E$ be any reduced component of $E$. Then either $\eta(E')=\{pt\}$ is a point (*) or $\eta(E')$ has codimension $2$ (**). 

By Corollary \ref{c:codim2} there exists a closed subset $Z\subset X$ with $\codim_X(Z)\geq 3$ such that every reflexive $p$-form $\s$ on $X\bsl Z$ extends to a regular $p$-form on $\eta^{-1}(X\bsl Z)$, for each $0\leq p\leq 3$. Thus for each divisor $E'$ of type (**) the following holds: Let $e\in E'$ be a generic point of $E'$, then $\eta^*\s$ is regular in $e\in E'$. Let $e'in E'$ be an arbitrary point of $E'$. Since $\X$ is smooth we can find an open neighbourhood $U\subset \X$ of $e'$ where $\Om^p_{\X}(\slog E)$ is free and $\eta^*\s$ is regular in a generic point of $U$. Then by \citep[Thm.~6.45]{GW10} $\eta^*\s$ is regular in $e'\in E'$. Since $e'\in E'$ was an arbitrary point, $\eta^*\s$ extends to all of $E'$ as a regular form.

It is only left to check that $\eta^*\s$ extends as a regular form on a divisor $E'$ of type (*). This case however is covered by Corollary \ref{c:12formsaffine}. \eproof

\section{Proof of Theorem \ref{t:2}}\label{ch. proofof1.2}
In this chapter we are going to prove Theorem \ref{t:2}. First let us recall the statement.

\begin{thmn}[Extension of $1$- and $2$-forms on good quotients]
Let $G$ be a reductive group and $V$ a smooth $G$-variety admitting a good quotient $X:=V\hq G$ of dimension $\dim(X)=n$. Let $\eta:\X\rightarrow X$ be a log resolution and $\s\in H^0\big(X,\,\Om_{X}^{[p]}\big)$ a reflexive $p$-form on $X$, for $p\in\{1,2\}$. Then 
$$
\eta^*(\s)\in H^0\big(\X,\,\Om_{\X}^{p}\big).
$$

\end{thmn}

To prove this Theorem we need the following lemma by Namikawa.

\begin{lem} \label{l:E't=0}
Let $X$ be a normal, affine variety with rational singularities and $x\in X$ a point in $X$. Let $\eta:\X\rightarrow X$ be a resolution of $X$ such that $E:=(\eta^{-1}(x))_{\mathrm{red}}$ is an snc divisor on $\X$. Then $H^0\big(E,\,\hat{\Om}_{E}^{p}\big)=0$ for all $p>0$.
\end{lem}

\proof
We outline a proof that uses a Hodge-theoretic method. More details on this argument can be found in \citep[Lemma~1.2]{Nam01a}. Recall that the cohomology group $H^{i}\big(E,\,\CC\big)$ carries a mixed Hodge structure. Assume that $H^0\big(E,\,\hat{\Om}_{E}^{p}\big)\neq 0$ for some $p>0$. Then by Hodge symmetry we get that $H^p\big(E,\,\hat{\Om}_{E}^{0}\big)=H^p\big(E,\,\mathcal{O}_{E}\big)\neq 0$. However, since $X$ has rational singularities $H^p\big(\X,\,\mathcal{O}_{\X}\big)=0$, from which one can prove that $H^p\big(E,\,\mathcal{O}_{E}\big)= 0$. This contradicts the assumption. \eproof

\proofof{Theorem \ref{t:2}}
The organisation of the following proof follows \citep[pp.~10-12]{Nam01}. Without loss of generality $V$ is affine and $X$ is a an affine GIT-quotient (see Lemma \ref{l:redtoaff}). By $E:=\exc(\eta)$ we denote the exceptional divisor. Recall that $E$ is an snc divisor.\\

\noindent\textbf{Step 1 (Preparations):}  
Let $n:=\dim(X)$ be the dimension of $X$. We denote by $X_{\mathrm{sing}}$ the singular locus of $X$. By Corollary \ref{c:codim2} the extension of $\eta^*(\s)$ is clear outside a certain locus $T\subset X_{\mathrm{sing}}$ of $codim_{X}(T)\geq 3$. Outside $T$ the GIT-quotient $X$ has finite quotient singularities only (see Proposition \ref{p:codim2}).  Thus it is left to check the extension of $\eta^*(\s)$ over $T$. Let $E^0\subset E$ be an irreducible component of $E$ with $T^0:=\eta(E^0)\subset T$. We put $k:=\dim(T)- \dim (T^0)$ and prove the extension of $\eta^*(\s)$ along $E^0$ by induction on $k$.  \\

\noindent\textbf{Step 2 (Case $k=0$):}
In this case $\dim(T)= \dim(T^0)$. We set $l:=\codim_{X}(T)$. Recall that $l\geq 3$. \\

\noindent\textbf{Step 2.1 (Cutting down):}
As in Lemma \ref{l:cutting down} we consider $n-l$ general hyperplanes $H_{1},$ $\dots,$ $H_{n-l}$ and denote by $H:=H_{1}\cap\dots\cap H_{n-l}$ a general $l$-dimensional complete intersection in $X$. Let $t_0\in T^0\cap H$. Using Lemma \ref{l:cutting down} we know the following facts. By replacing $X$ by a small open neighbourhood of $t_0$ we may assume that $T^0\cap H=\{t_0\}$. The preimage $\widetilde{H}:=\eta^{-1}(H)$ is a resolution of singularities. Since $X$ is an affine GIT-quotient, $H=\pi^{-1}(H)\hq G$ is an affine GIT-quotient as well. Thus $H$ has a unique distinguished singular point $t_0$ such that $H_{\mathrm{sing}}\bsl\{t_0\}$ contains only finite quotient singularities. Similar to the proof of Proposition \ref{p:codim2} we can replace $X$ by a smaller neighbourhood of $t_0$ such that there exists a variety $S$ and a flat surjective morphism $\F:X\rightarrow S$ with connected fibers, such that the restriction $\F_{\mid T^0}:T^0\rightarrow S$ is an isomorphism and such that there exist a point $s_0\in S$ with $\F^{-1}(s_0)=H$. By choosing $X$ small enough we can also ensure that $\Om_s^p$ is trivial for all $0\leq p\leq \dim(S)$. Given any point $s\in S$ the preimage $\Xs:=\F^{-1}(s)$ is an $l$-dimension affine GIT-quotient $X_s=\pi^{-1}(X_s)\hq G$. If $t_s\in T^0\cap X_s$ is the unique point that maps to $s\in S$, then $\Xs$ has finite quotient singularities away from $t_s\in \Xs$. The map $\eta:\X\rightarrow X$ gives a simultaneous resolution of the fibers $X_s$ for $s\in S$. Let $E'\subset \X$ be the union of all irreducible components of $E$ that map into $T^0$. Then $E'\rightarrow S$ is a proper map and by Lemma \ref{l:cutting down} every fiber $E'_s$ is a simple normal crossing divisor. Note that $E'_s$ has support $\eta^{-1}(t_s)$, for all $s\in S$.\\

\noindent\textbf{Step 2.2 (Filtrations):}
We consider the composition $\varPsi:= \F\circ\eta:\X\rightarrow S$. After shrinking $S$ we may assume, that $\varPsi$ is smooth. By \citep[Prop.~3.11]{Keb13} there are filtrations
$$
\hat{\Om}^q_{E'}=\mathcal{F}^0\supset\mathcal{F}^1\supset\dots\supset\mathcal{F}^q\supset\mathcal{F}^{q+1}=\{0\}
$$
and
$$
\Om^q_{\X}(\slog E')=\mathcal{G}^0\supset\mathcal{G}^1\supset\dots\supset\mathcal{G}^q\supset\mathcal{G}^{q+1}=\{0\}
$$
for $0\leq q\leq \dim(\X)-1$, where $\hat{\Om}^q_{E'}$ is the sheaf of torsion free $q$-forms on $E'$. These filtrations induce exact sequences

$$
0\rightarrow \mathcal{F}^{r+1}\rightarrow \mathcal{F}^r \rightarrow (\varPsi_{\mid E'})^*\Om_{S}^{r}\otimes \hat{\Om}^{q-r}_{E'/S}\rightarrow 0
$$
on $E'$ and
$$
0\rightarrow \mathcal{G}^{r+1}\rightarrow \mathcal{G}^r \rightarrow \varPsi^*\Om_{S}^{r}\otimes \Om^{q-r}_{\X/S}(\slog E')\rightarrow 0
$$
on $\X$ for $0\leq r\leq q$, where $\hat{\Om}^{q-r}_{E'/S}$ and $\Om^{q-r}_{\X/S}$ are sheaves of relative torsion free $(q-r)$-forms over $S$ (see \citep[Ch.~2.1]{Keb13}). In this proof we only have to consider $q\in\{1,\,2\}$. For $q=1$, we have $\mathcal{F}^{1}\cong (\varPsi_{\mid E'})^*\Om_{S}^{1}$ and $\mathcal{G}^{1}\cong \varPsi^*\Om_{S}^{1}$ and get one exact sequence each:

\begin{align}
\begin{split}
 & 0\rightarrow (\varPsi_{\mid E'})^*\Om_{S}^{1}\rightarrow \hat{\Om}^1_{E'}\rightarrow \hat{\Om}^{1}_{E'/S}\rightarrow 0\\
 & 0\rightarrow \varPsi^*\Om_{S}^{1}\rightarrow \Om^1_{\X}(\slog E') \rightarrow \Om^{1}_{\X/S}(\slog E')\rightarrow 0
\end{split}
\end{align}

For $q=2$, we have $\mathcal{F}^{2}\cong (\varPsi_{\mid E'})^*\Om_{S}^{2}$ and $\mathcal{G}^{2}\cong \varPsi^*\Om_{S}^{2}$. We set $\mathcal{F}:=\mathcal{F}^{1}$ and $\mathcal{G}:=\mathcal{G}^{1}$ and get two exact sequences each:
\begin{align}
\begin{split}
 & 0\rightarrow \mathcal{F} \rightarrow \hat{\Om}^2_{E'}\rightarrow \hat{\Om}^{2}_{E'/S}\rightarrow 0\\
 & 0\rightarrow \mathcal{G}\rightarrow \Om^2_{\X}(\slog E') \rightarrow \Om^{2}_{\X/S}(\slog E')\rightarrow 0
\end{split}
\end{align}
and 
\begin{align}
\begin{split}
 & 0\rightarrow (\varPsi_{\mid E'})^*\Om_{S}^{2}\rightarrow \mathcal{F}\rightarrow (\varPsi_{\mid E'})^*\Om_{S}^{1}\otimes \hat{\Om}^{1}_{E'/S}\rightarrow 0\\
 & 0\rightarrow \varPsi^*\Om_{S}^{2}\rightarrow \mathcal{G} \rightarrow \varPsi^*\Om_{S}^{1}\otimes \Om^{1}_{\X/S}(\slog E')\rightarrow 0
\end{split}
\end{align}

\noindent\textbf{Step 2.3 (The long exact cohomology sequences):}
Let us now consider the sequence (+), that we have already seen in Chapter \ref{sequence}
$$
0\rightarrow \Om_{\X}^{p}/ \Om_{\X}^{p}\big( \log E' \big)\big( - E' \big)\rightarrow \Om_{\X}^{p}\big( \log E' \big)/ \Om_{\X}^{p}\big( \log E' \big)\big( - E' \big) \rightarrow \Om_{\X}^{p}\big( \log E' \big)/ \Om_{\X}^{p}\rightarrow 0.
$$
We want to prove the following claim.\\

\noindent\textbf{Claim 1.}
The induced maps
$$
\gamma_p: H^0\big(E',\,\Om_{\X}^{p}/ \Om_{\X}^{p}\big( \slog E' \big)\big( - E' \big)\big)\rightarrow H^0\big(E',\,\Om_{\X}^{p}\big( \slog E' \big)/ \Om_{\X}^{p}\big( \slog E' \big)\big( - E' \big) \big)
$$
are surjective for $p\in\{1,\, 2\}$.\\

\noindent\proofof{Claim 1}
Recall from Chapter \ref{Ch.MHS} that $\Om_{\X}^{p}/ \Om_{\X}^{p}\big( \slog\, E' \big)\big( - E' \big)\cong \hat{\Om}^{p}_{E'}$ for all $0\leq p \leq \dim(\X)$. By tensoring the exact sequence 
$$
0\rightarrow \mathcal{O}_{\X}(- E')\rightarrow \mathcal{O}_{\X} \rightarrow \mathcal{O}_{E'}\rightarrow 0
$$
with $\Om_{\X}^p(\slog E')$, we see that $\Om_{\X}^{p}\big(\slog E' \big)/ \Om_{\X}^{p}\big(\slog E' \big)\big( - E' \big)\cong \Om^{p}_{\X}(\slog E')\otimes\mathcal{O}_{E'}$ for all $0\leq p \leq \dim(\X)$. This induces the natural map 
$$ \hat{\Om}^{p}_{E'}\cong\Om_{\X}^{p}/ \Om_{\X}^{p}\big(\slog E' \big)\big( - E' \big)\rightarrow \Om_{\X}^{p}\big( \slog E' \big)/ \Om_{\X}^{p}\big(\slog  E' \big)\big( - E' \big)\cong\Om^{p}_{\X}(\slog E')\otimes\mathcal{O}_{E'},$$
for all $0\leq p\leq \dim(\X)$. Similarly we get a natural map between the sheaves 
$$
\hat{\Om}^{p}_{E'/S}\rightarrow \Om^{p}_{\X/S}(\slog E')\otimes\mathcal{O}_{E'}.
$$
Thus, if we tensor the second sequence in ($1$) with $\mathcal O_{E'}$ we have a map from each element in the first sequence in ($1$) to the second sequence, which induces the following commutative diagram with exact columns:
\begin{align}\tag{I}
\begin{split}
\xymatrix{
0\ar[d] & 0\ar[d] \\
H^0\big(S,\,\Om_{S}^{1} \big) \ar[d]\ar[r]  & H^0\big(S,\,\Om_{S}^{1} \big)\ar[d]\\
H^0\big(E',\,\hat{\Om}^1_{E'}\big)\ar[d]\ar[r]^{\gamma_1~~~~~~~~} & H^0\big(E',\,\Om^1_{\X}(\slog E')\otimes\mathcal{O}_{E'}\big)\ar[d]\\
H^0\big(E',\,\hat{\Om}^{1}_{E'/S}\big) \ar[d]\ar[r]^{\mu_1~~~~~~~~} & H^0\big(E',\,\Om^{1}_{\X/S}(\slog E')\otimes\mathcal{O}_{E'}\big)\ar[d]\\
\dots & \dots.
}
\end{split}
\end{align}
In this diagram we can identify $H^0\big(E',\,\varPsi^*\Om_{S}^{1}\otimes \mathcal O_{E'} \big)\cong H^0\big(S,\,\Om_{S}^{1} \big)$, because $\varPsi^*_{\mid E'}:E'\rightarrow S$ is a proper map.

Analogous, tensoring the second sequence in each ($2$) and ($3$) with $\mathcal O_{E'}$ we get the following two commutative diagrams with exact columns:
\begin{align}\tag{II}
\begin{split}
\xymatrix{
0\ar[d] & 0\ar[d] \\
H^0\big(E',\,\mathcal F \big) \ar[d]\ar[r]  & H^0\big(E',\,\mathcal G \otimes\mathcal{O}_{E'} \big)\ar[d]\\
H^0\big(E',\,\hat{\Om}^2_{E'}\big)\ar[d]\ar[r]^{\gamma_2~~~~~~~~} & H^0\big(E',\,\Om^2_{\X}(\slog E')\otimes\mathcal{O}_{E'}\big)\ar[d]\\
H^0\big(E',\,\hat{\Om}^{2}_{E'/S}\big) \ar[d]\ar[r]^{\mu_2~~~~~~~~} & H^0\big(E',\,\Om^{2}_{\X/S}(\slog E')\otimes\mathcal{O}_{E'}\big)\ar[d]\\
\dots & \dots.
}
\end{split}\\
\nonumber
\end{align}
\begin{align}\tag{III}
\begin{split}
\xymatrix{
0\ar[d] & 0\ar[d] \\
H^0\big(S,\,\Om_{S}^{2} \big) \ar[d]\ar[r]  & H^0\big(S,\,\Om_{S}^{2} \big)\ar[d]\\
H^0\big(E',\,\mathcal F \big) \ar[d]\ar[r]  & H^0\big(E',\,\mathcal G \otimes\mathcal{O}_{E'} \big)\ar[d]\\
H^0\big(E',\,\hat{\Om}^{1}_{E'/S}\big)^{n-l} \ar[d]\ar[r]^{(\mu_1)^{n-l}~~~~~~~~} & H^0\big(E',\,\Om^{1}_{\X/S}(\slog E')\otimes\mathcal{O}_{E'}\big)^{n-l}\ar[d]\\
\dots & \dots.
}
\end{split}
\end{align}
In diagram (III) we were allowed to exchange the map 
$$
H^0\big(E',\,(\varPsi_{\mid E'})^*\Om_{S}^{1}\otimes\hat{\Om}^{1}_{E'/S}\big)\rightarrow H^0\big(E',\,(\varPsi^*\Om_{S}^{1}\otimes\hat{\Om}^{1}_{\X/S}\otimes \mathcal O_{E'}\big)
$$
by $(\mu_1)^{n-l}$, because $\Om_{S}^{1}$ is trivial. Here, $n-l$ is the dimension of $S$.

We first want to show that $H^0\big(E',\,\hat{\Om}^{p}_{E'/S}\big) = H^0\big(E',\,\Om^{p}_{\X/S}(\slog E')\otimes\mathcal{O}_{E'}\big)=0$ for $p\in\{1,\,2\}$. Then the assertion of the Claim $1$ for $p=1$ follows from (I) and for $p=2$ from (II) and (III).\\

\noindent\textbf{Claim 2.}
$H^0\big(E',\,\hat{\Om}^{p}_{E'/S}\big)=0$ for $p\in \{1,\, 2\}$.\\

\noindent\proofof{Claim 2}

By Lemma \ref{l:E't=0} we know that $H^0\big(E'_{s},\,\hat{\Om}^{p}_{E'_{s}}\big)=0$ for $s\in S$ and $p\in \{1,\, 2\}$. We assume that there exists a non-zero section $\tau\in H^0\big(E',\,\hat{\Om}^{p}_{E'/S}\big)$ and use the assertion of Lemma \ref{l:E't=0} to get a contradiction. If $s\in S$ is general then the restriction $\tau_{\mid E'_s}$ does not vanish, 
$$ \tau_{\mid E'_s} \in H^0\big(E'_s,\,(\hat{\Om}^{p}_{E'/S})_{\mid E'_s}\big)\bsl\{0\}.$$
By \citep[Cor.~3.10]{Keb13} we have the following isomorphism:
$$H^0\big(E'_s,\,(\hat{\Om}^{p}_{E'/S})_{\mid E'_s}\big)\cong H^0\big(E'_{s},\,\hat{\Om}^{p}_{E'_{s}}\big).$$ 
However, by Lemma \ref{l:E't=0} the right-hand side is zero. This contradicts the assumption and ends the proof of Claim $2$.\hfill Q.E.D.\\

\noindent\textbf{Claim 3.}
$H^0\big(E',\, \Om^{p}_{\X/S}(\slog E')\otimes\mathcal{O}_{E'}\big)=0$ for all $p\in \{1,\, 2\}$.\\

\noindent\proofof{Claim 3}
Recall from the remark to Lemma \ref{c:12forms} that the maps $\mu_p(s):H^0\big(E'_{s},\,\hat{\Om}^p_{E'_{s}}\big)\rightarrow H^0\big(E'_{s},\,\Om^p_{\X_{s}}(\slog E'_{s})\otimes\mathcal{O}_{E'_{s}}\big)$ are surjective for all $s\in S$ and $p\in \{1,\, 2\}$. Thus $H^0\big(E'_{s},\,\hat{\Om}^{p}_{E'_{s}}\big)=0$ implies that $H^0\big(E'_{s},\,\Om^{p}_{\X_{s}}(\slog E'_{s})\otimes\mathcal{O}_{E'_{s}}\big)=0$, hence $H^0\big(E',\, \Om^{p}_{\X/S}(\slog E')\otimes\mathcal{O}_{E'}\big)=0$ for all $p\in \{1,\, 2\}$ by the same argument as in the proof of Claim $2$. \hfill Q.E.D.\\

Now we can finish the proof of Claim $1$ by analysing the cases $p=1$ and $p=2$ separately. \\

\noindent\underline{Case $p=1$:}
Consider the diagram (I). The first horizontal map $H^0\big(S,\Om^1_S\big)\rightarrow H^0\big(S,\Om^1_S\big)$ is surjective. Thus Claim $2$ and $3$ ensure that $\gamma_1$ is surjective as well. \\

\noindent\underline{Case $p=2$:}
First consider diagram (III). Since the top map $H^0\big(S,\Om^2_S\big)\rightarrow H^0\big(S,\Om^2_S\big)$ is surjective, Claim $2$ and $3$ ensure that $H^0\big(E',\,\mathcal F \big) \rightarrow H^0\big(E',\,\mathcal G \otimes\mathcal{O}_{E'} \big)$ is surjective as well. This map is the first horizontal map in diagram (II). By the same argument as in the case $p=1$ using Claim $2$ and $3$ we see that $\gamma_2$ is surjective. This finishes the proof of Claim $1$.\hfill Q.E.D.\\
 
\textbf{Step 2.4 (Proof of case $k=0$):}
By taking the cohomology of the exact sequence 
$$
0\rightarrow \Om_{\X}^{p}/ \Om_{\X}^{p}\big(\log E' \big)\big( - E' \big)\rightarrow \Om_{\X}^{p}\big(\log E' \big)/ \Om_{\X}^{p}\big(\log E' \big)\big( - E' \big) \rightarrow \Om_{\X}^{p}\big(\log E' \big)/ \Om_{\X}^{p}\rightarrow 0
$$
and by applying Claim $1$ of Step $2.3$ we see that the boundary map 
$$
\delta: H^0\big(E',\,\Om_{\X}^{p}\big(\log E' \big)/ \Om_{\X}^{p} \big) \rightarrow H^1\big(E',\,\Om_{\X}^{p}/\Om_{\X}^{p} \big(\log E' \big)\big( - E' \big) \big)
$$
in the long exact sequence is an injection for $p\in \{1,\, 2\}$. Similar to the proof of Lemma \ref{deltainjective} we get that $\tau: H^0\big(\X,\,\Om_{\X}^{p} \big) \rightarrow H^0\big(\X,\,\Om_{\X}^{p}\big(\log E' \big)\big)$ is surjective. This finishes the proof of the case $k=0$. \\

\noindent\textbf{Step 3 (Case $k>0$):}
We assume that the result holds for all cases where $0\leq k\leq m$ for some $m\in\NN$ and want to show it for $k=m+1$.  We once again set $l:=\codim_{X}(T)$. \\

\noindent\textbf{Step 3.1 (Cutting down):}
As in Step $2.1$ we consider $n-l-k$ general hyperplanes $H_{1},$ $\dots,$ $H_{n-l-k}$ and denote by $H:=H_{1}\cap\dots\cap H_{n-l-k}$ a general $(l+k)$-dimensional complete intersection in $X$. Let $t_0\in T^0\cap H$. As before we know the following facts. By replacing $X$ by a small open neighbourhood of $t_0$ we may assume that $T^0\cap H=\{t_0\}$. The preimage $\widetilde{H}:=\eta^{-1}(H)$ is a resolution of singularities. Since $X$ is an affine GIT-quotient $H$ is an affine GIT-quotient as well. As in Step $2.1$ we can define a flat surjective morphism $\F:X\rightarrow S$ with connected fibers, such that the restriction $\F_{\mid T^0}:T^0\rightarrow S$ is an isomorphism and such that there exist a point $s_0\in S$ with $\F^{-1}(s_0)=H$. Given any point $s\in S$ the preimage $\Xs:=\F^{-1}(s)$ is a $(l+k)$-dimension affine GIT-quotient $X_s$. By $t_s\in T^0\cap X_s$ we denote the unique point that maps to $s\in S$. The map $\eta:\X\rightarrow X$ gives a simultaneous resolution of the fibers $X_s$, for $s\in S$. Let $E'\subset \X$ be the union of all irreducible components of $E$ that map into $T^0$. Then $E'\rightarrow S$ is a proper map and by Lemma \ref{l:cutting down} every fiber $E'_s$ is a simple normal crossing divisor. Note that $E'_s$ has support $\eta^{-1}(t_s)$, for all $s\in S$. \\

\noindent\textbf{Step 3.2 (Proof of case $k>0$):}
By the induction hypothesis $H^0\big(\X,\,\Om_{\X}^{p}\big(\log E \big)\big) \cong H^0\big(\X,\,\Om_{\X}^{p}\big(\log E' \big)\big)$ for all $p\in\{1,\,2\}$. Thus we only have to prove that  $\tau: H^0\big(\X,\,\Om_{\X}^{p} \big) \rightarrow H^0\big(\X,\,\Om_{\X}^{p}\big(\log E' \big)\big)$ is surjective for $p\in\{1,\,2\}$. Following the argumentation of Step $2.2$ to Step $2.4$ we get the result. 
\eproof

\part{The 4-dimensional case}\label{c: dim=4}

In this part we want to analyse Theorem \ref{t:1} when $X$ is a $4$-dimensional good quotient of a smooth variety by a reductive group. Let $\eta:\X\rightarrow X$ be a log resolution and $\s\in H^0\big(X,\,\Om_{X}^{[p]}\big)$ a reflexive $p$-form on $X$, $0\leq p\leq 4$. For $p=0$ and $p=4$, the extension of $\eta^*(\s)$ to all of $\X$ without poles follows from the fact that $X$ has rational singularities (see Remark to Proposition \ref{nforms}). For $p\in\{1,\,2\}$, the result was proven in the last chapter. Thus, it is only left to check the Extension Theorem for $3$-forms. 

\begin{rem}
Since $\dim(X)=n=4$, in this part we are interested in the special case of an $(n-1)$-form. When $(X,D)$ is a klt-pair this case is handled separately (see \citep[Prop.~6.1]{GKK10}). However, those kind of arguments won't work in our situation, because of two reasons. First of all we do not know anything about the discrepancy (see \citep[Prop.~5.1]{GKK10}). Secondly we do not have a duality between $(n-1)$-forms and vector fields (see \citep[Prop.~6.1]{GKK10}), because in general $K_X$ is not Cartier and thus $\mathcal O_X(K_X)\not\cong\mathcal O_X$.
\end{rem}

By Proposition \ref{logpoles} we already know that $\eta^*(\s)$ extends with log-poles along the exceptional divisor $E:= \exc(\eta)$. In this case the residue sequence introduced in \citep{EV92} is a useful tool to show that $\eta^*(\s)$ in fact extends without these poles. The details of this idea will be explained in Chapter \ref{subs:residuemap}. Examples how this sequence is used in the case of a klt-pair $(X,D)$ can be found in \citep[Rem.~11.8, 25.E]{GKKP11}. Similar to the proof in the klt case, we would like to split the log resolution $\eta:\X\rightarrow X$ into a finite sequence of surjective, birational morphisms
\begin{align}\tag{$\star$}
\X\rightarrow X_k\rightarrow \dots \rightarrow X_1 \rightarrow X_0:=X,
\end{align}
such the following properties hold.
\enumeraten
\begin{enumerate}
\item The morphism $\X\rightarrow X_k$ is a log resolution and $X_k$ has finite quotient singularities.
\item Each morphism $X_{i+1} \rightarrow X_{i}$ contracts exactly one irreducible divisor that is a strict transform of an irreducible component of $E$ that is not contracted by $\X\rightarrow X_k$.
\end{enumerate}
Then we can apply a modified version of the residue sequence to each step separately. The desired sequence ($\star$) is induced by the partial resolution algorithm of Kirwan and will be introduced in the next chapter. A modified version of the residue sequence that can be applied in this case will be introduced in Chapter \ref{subs:residuemap}.

\section{The partial resolution of Kirwan}\label{ch. Kirwan}
Let $V$ be a smooth projective variety and suppose that a reductive group $G$ acts linearly on $V$ with respect to a projective embedding $V\subset \PP^N$. Then $V^{ss}\hq G$ is a projective GIT-quotient in the sense of Definition \ref{d:projGIT} (i.e. the $G$-linearisation $\mathcal L$ on $V$, corresponding to the projective embedding, is ample). Throughout this chapter, whenever we consider a $G$-action on a projective variety $V$, we think of an action of this type. By abuse of notation we will write $X=V\hq G$ instead of $X=V^{ss}\hq G$ to simplify the diagrams in the following chapters.

\begin{thm}[Kirwan's partial Resolution, {\citep[6.1, 6.3, 6.9]{Kir85}}]\label{T:KirwanRes}
Let $G$ be a reductive group and $V$ a smooth projective $G$-variety with a $G$-linearisation $\mathcal L=\mathcal L_0$ on $V$, such that we get a projective GIT-quotient $X:=V\hq G$. If $V^s\neq \emptyset$, then there exists a finite sequence of blow-ups
\begin{align}\tag{K1}
\xymatrix{
V_k \ar[r]_{F_k}\ar@/^1pc/[rrrr]^{F}  & V_{k-1} \ar[r]_{F_{k-1}} & \dots \ar[r]_{F_{2}} & V_1 \ar[r]_{F_{1}} & V_0=V , & \text{for }k\in\mathbb{N},
}
\end{align}
along smooth $G$-invariant subspaces $W_i\subset V_i$, for $0\leq i \leq k-1$, such that we have a $G$-action on each $V_{i+1}$ (lifting the action on $V_i$) and a $G$-linearisation $\mathcal L_{i+1}$ on $V_{i+1}$ with the following properties:
\begin{enumerate}
\item Every semi-stable point in $V_k$ is stable (with respect to $\mathcal L_{k}$).
\item For all $0\leq i\leq k-1$, every lift in $V_{i+1}$ of a stable (respectively not semi-stable) point in $V_i$ (with respect to $\mathcal L_{i}$) is stable (respectively not semi-stable) in $V_{i+1}$ (with respect to $\mathcal L_{i+1}$).
\end{enumerate} 

Moreover, this sequence induces a sequence of blow-ups of the induced GIT-quotients
\begin{align}\tag{K2}
\xymatrix{
X_k \ar[r]^{f_k}\ar@/_1pc/[rrrr]_{f}  & X_{k-1} \ar[r]^{f_{k-1}} & \dots \ar[r]^{f_{2}} & X_1 \ar[r]^{f_{1}} & X_0=X 
}
\end{align}
along the images $W_i\hq G$ of $W_i$ in $X_i:=V_i\hq G$, for $0\leq i\leq k-1$, such that the diagram 
\begin{align}\tag{K3}
\xymatrix{
**[l]V_k^{ss} \ar[r]_{F_k^{ss}}\ar@/^1pc/[rrrr]^{F^{ss}}\ar[d]_{\pi_{k}}  & V_{k-1}^{ss} \ar[r]_{F_{k-1}^{ss}}\ar[d]_{\pi_{k-1}} & \dots \ar[r]_{F_{2}^{ss}} & V_1^{ss} \ar[r]_{F_{1}^{ss}}\ar[d]_{\pi_{1}} & **[r] V_0^{ss}:=V^{ss}\ar[d]_{\pi_0}^{=~~\pi}\\
**[l]X_k \ar[r]^{f_k}\ar@/_1pc/[rrrr]_{f}  & X_{k-1} \ar[r]^{f_{k-1}} & \dots \ar[r]^{f_{2}} & X_1 \ar[r]^{f_{1}} & **[r]  X_0:=X 
}
\end{align}
commutes. Here, $F_i^{ss}$ and $F^{ss}$ denote the restrictions of $F_i$ and $F$ to the semi-stable locus. We will refer to this as the \emph{partial resolution of Kirwan}.

\end{thm}

\begin{rem}
The following remaks should help the reader to understand the diagram (K3) better and should illustrate the fact that we can construct the sequence ($\star$), presented in the introduction to Part \ref{c: dim=4} from (K$3$).

\enumeraten
\begin{enumerate}
   \item In this paper we want $k\in \ZZ$ to be the minimal number such that the assertion of the previous theorem are true. We call $k$ the \emph{number of Kirwan steps}.
   \item Let $0\leq i\leq k-1$ be an integer. The $G$-invariant and smooth subset $W_i\subset V_i$ contains all semi-stable points in $V_i$ (with respect to $\mathcal L_i$) that have a stabiliser isomorphic to a reductive group $R_i\subset G$ of maximal dimension (see \citep[Ch.~5]{Kir85}). By \citep[Lem.~6.1]{Kir85} we know that every semi-stable point in the $F_{i+1}$-exceptional locus (with respect to $\mathcal{L}_{i+1}$) is no longer fixed by a conjugate of $R_i$.
\item Since $V_k^{ss}=V_k^{s}$ the quotient $\pi_k:V_k^{ss}\rightarrow X_k$ is a geometric quotient and $X_k$ has only finite quotient singularities (see Example \ref{e:fqs stable}). Since $f:X_k\rightarrow X$ is surjective we can call $f$ a partial resolution.
 \item Let $E_f\subset X_k$ be the $f$-exceptional divisor in $X_k$. Then $E_f$ consists of $k$ irreducible components $E_1$, $\dots$, $E_k$, where $E_i$ is the strict transform of the $f_i$-exceptional divisor, for $1\leq i\leq k$.

\end{enumerate}
\end{rem}

\begin{lem}[The exceptional locus and fibers of the blow-ups on the level of quotients] \label{l:fiber of Kirwan Step}
In the setting of Theorem \ref{T:KirwanRes}, let $1\leq i\leq k$ be a positive number. Consider the $i$-th Kirwan step with the same notation as in the previous Theorem and the corresponding diagram:
$$
   \xymatrix{
    V_i^{ss} \ar[r]^{F_{i}^{ss}}\ar[d]_{\pi_{i}} & **[r] V_{i-1}^{ss}\ar[d]^{\pi_{i-1}}\\
     X_i \ar[r]_{f_{i}} & **[r]  X_{i-1}. 
   }
$$
\enumeratel

\begin{enumerate}
\item 
Let $T_{i-1}$ be the image $W_{i-1}\hq G$ of $W_{i-1}$ in $X_{i-1}$. Then the exceptional locus $f_i^{-1}(T_{i-1})$ of $f_i$ is a quotient $F^{-1}_i(W_{i-1})\hq G$ (with respect to the restriction of the linearisation $\mathcal L_i$ to $F_i^{-1}(W_{i-1})$) of a smooth projective variety $F_i^{-1}(W_{i-1})$ by $G$.  
\item
Let $x_{i-1}\in T_{i-1}$ be a point in the quotient. Then there exists a $m\in \mathbb N_+$ such that $f_i^{-1}(x_{i-1})$ is isomorphic to a projective GIT-quotient $\PP^m\hq H$ (with respect to some $H$-linearisation on $\PP^m$), where $H\subset G$ is a reductive subgroup of $G$.
\end{enumerate} 

\end{lem}

\proof
The first assertion follows from the definition of a blow-up of a smooth variety along a smooth subspace and \citep[Lem.~3.11, Rem.~6.8]{Kir85}. Since $W_{i-1}$ is smooth and $G$-invariant the preimage $F_{i}^{-1}(W_{i-1})$ is smooth with an induced $G$-action. By the construction of diagram ($K3$) there exists a $G$-linearisation on $F_{i}^{-1}(W_{i-1})$, such that the GIT-quotient $F_i^{-1}(W_{i-1})\hq G$ coincides with the $f_{i}$-exceptional divisor (see \citep[Rem.~6.8]{Kir85}).

The proof of the second assertion is based on \citep[Rem.~6.4]{Kir85}. Let $x_{i-1}\in T_{i-1}$ be a point in the quotient. Then there exists a $v_{i-1}\in\pi_{i-1}^{-1}(x_{i-1})\subset W_{i-1}$ with closed $G$-orbit. Thus, the preimage $F_{i}^{-1}(v_{i-1})$ is isomorphic to the projectification $\mathbb{P}(N_{v_{i-1}})$ of the normal space $N_{v_{i-1}}$ to $W_{i-1}$ at the point $v_{i-1}$. The stabiliser $G_{v_{i-1}}$ acts linearly on $\mathbb{P}(N_{v_{i-1}})$. This action is induced by its canonical action of on the normal space $N_{v_{i-1}}$. Since $v_{i-1}\in W_{i-1}\cap V_{i-1}^{ss}$ the group $G_{v_{i-1}}$ is a reductive subgroup of $G$ (see previous remark, assertion $2$). As before there exists a $G_{v_{i-1}}$-linearisation on $\mathbb{P}(N_{v_{i-1}})$, such that $\mathbb{P}(N_{v_{i-1}})\hq G_{v_{i-1}}$ coincides with the preimage $f^{-1}_{i}(x_{i-1})$ (see \citep[Rem.~6.8]{Kir85}).
\eproof

\begin{rem}
As we have seen in the previous lemma, the $f_i$-exceptional locus $f_i^{-1}(T_{i-1})$ and the fiber of a point $x_{i-1}$ under the blow-up $f_i$ can be expressed as GIT-quotients of smooth varieties. Thus, the exceptional locus and these fibers are normal varieties and the notion of a reflexive differential form on them is well defined.
\end{rem}

The partial resolution was formulated and proven only in the case where $V$ is a smooth projective $G$-variety with $V^s\neq\emptyset$. By Corollary \ref{p:redtoV} we can prove Theorem \ref{t:3} by reducing to the case of an affine GIT-quotient $X:=V\hq G$ of a vector space $V$ by a group $G$ with linear action on $V$. In general this vector space might not fulfil the condition $V^s\neq\emptyset$. Thus, we need to find a partial resolution algorithm for this case as well. Since we are only interested in the GIT-quotient $X:=V\hq G$, we can change $G$ and $V$ such that the partial resolution algorithm of Kirwan induces a new algorithm in this case, too. As a preparation we need the following two results.

\begin{lem}\label{l:saturated}
Let $Z$ be a smooth variety with an action of a reductive group $G$. Let $U\subset Z$ be an open $G$-invariant subvariety admitting a good quotient $\pi:U\rightarrow U\hq G$. If $U\hq G$ is quasi-projective, then there exists a $G$-linearisation on $Z$ such that $U$ is a $G$-saturated subset of $Z^{ss}$. 
\end{lem}

\proof
The result follows from \citep[Thm.~on~p.1]{Hau03}. In his paper Hausen works on a normal variety $Z$ instead of a smooth variety and thus uses a Weil divisor $D$ on $Z$ for the $G$-linearisation instead of a line bundle $\mathcal L$. In Definition $1.1$ he then explains that in the smooth case, where $D$ is Cartier (and thus corresponds to a line bundle $\mathcal L$), his definition of a $G$-linearisation induced by $D$ coincides with the original definition of a $G$-linearisation induced by $\mathcal L$ in \cite{MFK94}.
\eproof

\begin{cor}\label{c:saturated}
Let $V$ be a vector space and $G$ a reductive group acting linearly on $V$. Let $X:=V\hq G$ be the induced GIT-quotient. Then there exists a smooth projective $G$-variety $\overline{V}$, together with a $G$-linearisation on $\overline{V}$, such that $V$ is a $G$-saturated subset of $\overline{V}^{ss}$ and the $G$-action on $\overline{V}$ is linear with respect to some projective embedding. 
\end{cor}

\proof
We choose $\overline{V}:=\PP(V\oplus \CC)$ with the induced $G$-action. Then the result follows from Lemma \ref{l:saturated}. The $G$-action is linear with respect to some projective embedding since any chosen linearisation $\mathcal L$ on $\PP(V\oplus \CC)$ is ample. 
\eproof

\begin{prop}[Kirwan's partial resolution for GIT-quotients of vector spaces]\label{p:Kirwanforaffine}
Let $G$ be a reductive group acting linearly on a vector space $V$ and consider the induced affine GIT-quotient $X:=V\hq G$. Let $\mathcal L=\mathcal L_0$ be the trivial $G$-linearisation on $V$, such that $V^s\neq \emptyset$. Then there exists a finite sequence of blow-ups
\begin{align}\tag{K1v}
\xymatrix{
V_k \ar[r]_{F_k}\ar@/^1pc/[rrrr]^{F}  & V_{k-1} \ar[r]_{F_{k-1}} & \dots \ar[r]_{F_{2}} & V_1 \ar[r]_{F_{1}} & V_0=V , & \text{for }k\in\mathbb{N},
}
\end{align}
along smooth $G$-invariant subspaces $W_i\subset V_i$, for $0\leq i\leq k-1$, such that we have a $G$-action on each $V_{i+1}$ (lifting the action on $V_{i}$) and a $G$-linearisation $\mathcal L_{i+1}$ on $V_{i+1}$ with the following properties. 
\begin{enumerate}
\item Every semi-stable point in $V_k$ is stable (with respect to $\mathcal L_{k}$).
\item For all $0\leq i\leq k-1$, every lift in $V_{i+1}$ of a stable (respectively not semi-stable) point in $V_i$ (with respect to $\mathcal L_{i}$) is stable (respectively not semi-stable) in $V_{i+1}$ (with respect to $\mathcal L_{i+1}$).
\end{enumerate} 

Moreover, this sequence induces a sequence of blow-ups of the induced GIT-quotients
\begin{align}\tag{K2v}
\xymatrix{
X_k \ar[r]^{f_k}\ar@/_1pc/[rrrr]_{f}  & X_{k-1} \ar[r]^{f_{k-1}} & \dots \ar[r]^{f_{2}} & X_1 \ar[r]^{f_{1}} & X_0=X 
}
\end{align}
along the images $W_i\hq G$ of $W_i$ in $X_i:=V_i\hq G$, for $0\leq i\leq k-1$, such that the diagram 
\begin{align}\tag{K3v}
\xymatrix{
**[l]V_k^{ss} \ar[r]_{F_k^{ss}}\ar@/^1pc/[rrrr]^{F^{ss}}\ar[d]_{\pi_{k}}  & V_{k-1}^{ss} \ar[r]_{F_{k-1}^{ss}}\ar[d]_{\pi_{k-1}} & \dots \ar[r]_{F_{2}^{ss}} & V_1^{ss} \ar[r]_{F_{1}^{ss}}\ar[d]_{\pi_{1}} & **[r] V_0^{ss}:=V^{ss}\ar[d]_{\pi_0}^{=~~\pi}\\
**[l]X_k \ar[r]^{f_k}\ar@/_1pc/[rrrr]_{f}  & X_{k-1} \ar[r]^{f_{k-1}} & \dots \ar[r]^{f_{2}} & X_1 \ar[r]^{f_{1}} & **[r]  X_0:=X 
}
\end{align}
commutes. Here, $F_i^{ss}$ and $F^{ss}$ denote the restrictions of $F_i$ and $F$ to the semi-stable locus. We will refer to this as the \emph{partial resolution of Kirwan for vector spaces}.
\end{prop}

\proof
By Corollary \ref{c:saturated} there exists a smooth projective $G$-variety $\overline V$ togehter with some $G$-linearisation such that $V$ is a $G$-saturated subset in $\overline{V}^{ss}$ and the $G$-action on $\overline{V}$ is linear with respect to some embedding. Since $V^s\neq\emptyset$ it follows that $\overline{V}^s\neq\emptyset$. Thus, Kirwan's partial resolution algorithm works for $\overline V$ and induces the diagrams $(K1)-(K3)$. If we consider $(K3)$ only over $V \subset \overline{V}^{ss}$ (or $X=V\hq G$ respectively) the restriction defines the diagrams $(K1v)-(K3v)$ with the needed properties.
\eproof

\begin{rem}
By construction, the partial resolution of Kirwan for vector spaces has the same properties as the original partial resolution of Kirwan. In particular, the remark to Theorem \ref{T:KirwanRes} and Lemma \ref{l:fiber of Kirwan Step} are fulfilled in this setting as well.
\end{rem}

\begin{lem}
Let $G'$ be a reductive group acting linearly on a vector space $V'$. Let $X:=V'\hq G'$ be the induced  GIT-quotient. Then there exists a reductive group $G$ acting linearly on a vector space $V$ such that $X\cong V\hq G$ and such that, using the trivial $G$-linearisation on $V$, we have $V^s\neq \emptyset$. 
\end{lem}

\proof
The existence of $V$ and $G$ follows from Proposition \ref{l:lunarichardson} together with the remark to this proposition. By Corollary \ref{c:lunarichardson} we know that $V^s\neq\emptyset$.
\eproof

The previous lemma implies the following result:

\begin{cor}\label{c:Kirw ex}
Let $X$ be a good quotient of a smooth variety by a reductive group. Then we can find an affine GIT-quotient $V\hq G$, where $G$ is a reductive group acting linearly on a vector space $V$ with the following properties:
\enumeraten
\begin{enumerate}
 \item The Extension Theorem for $V\hq G$ implies the Extension Theorem for $X$.
 \item The partial resolution algorithm presented in Proposition \ref{p:Kirwanforaffine} works for $V\hq G$.
\end{enumerate} 
\end{cor}

\proof
The first assertion follows from Lemma \ref{l:redtoV} and Corollary \ref{p:redtoV}. The second assertion is a direct consequence of the previous lemma. \eproof

The following example illustrates the partial resolution algorithm of Kirwan. 

\begin{ex} 
Consider the action of the group $G=\mathbb{C}^*$ of invertible complex numbers on the $2$-dimensional complex space $V=\mathbb{C}^2$ via $t.(z_1,z_2):=(t\cdot z_1,t^{-1}\cdot z_2)$. We can write $V=\spec(A)$, where $A:=\CC[Z_1,Z_2]$ is the polynomial ring of $V$. The ring of $G$ invariant polynomial is $A^G=\CC[Z_1\cdot Z_2]$. Thus the quotient $X=V\hq G$ is isomorphic to the space of complex numbers $\mathbb{C}$.

If we consider the trivial linearisation on $V$ we get $V^{ss}=V$. One can check, that the axes of $\CC^2$ contain all the vectors $v\in V$ with $v\notin V^s$ and that $0\in V^{ss}$ is the semi-stable vector that has the largest stabilizer (with respect to the dimension). Since the orbit $G(0)$ is closed the space along which we have to blow up according to the algorithm of Kirwan is $W:=\{0\}$. 

Consider the blow-up $F: \widetilde V\rightarrow V$ of $V$ in $0\in V$. The fiber $F^{-1}(\{0\})$ is isomorphic to $\PP^1$. One can easily see that the points $v_1$, $v_2$ corresponding to $[1:0]$ and $[0:1]\in\PP^1$ are the only non-semi-stable points of $\widetilde V$, because polynomials of the form $Z_1\cdot Z_2$ vanish only in these points. Except for $v_1$ and $v_2$ every other point in $v\in\widetilde V$ has trivial stabiliser. Since $(\widetilde V)^{ss}=\widetilde V \bsl \{v_1,v_2\}$ the orbit $G(v)$ is closed in $(\widetilde V)^{ss}$. Thus $(\widetilde V)^{ss}=(\widetilde V)^{s}$ and the algorithm stops. The quotient $\widetilde V\hq G$ is isomorphic to $\CC$ and the induced blow-up of quotients $f:\widetilde V\hq G\rightarrow V\hq G$ is an isomorphism. 
\end{ex}

\section{Preparation for Theorem \ref{t:3}}\label{Ch.prep for 3}
In this chapter we want to present some reduction steps and results that will help us proving Theorem \ref{t:3}. For the convenience of the reader lets recall the theorem.

\begin{thmn}[Extension of $(n-1)$-forms on good quotients]
Let $G$ be a reductive group and $V$ a smooth $G$-variety admitting a good quotient $X:=V\hq G$ of dimension $n\geq 1$. Let $\eta:\X\rightarrow X$ be a log resolution and $\s\in H^0\big(X,\,\Om_{X}^{[n-1]}\big)$ a reflexive $(n-1)$-form on $X$. Assume that the Extension Theorem is true for all reflexive $p$-forms on any good quotient of dimension less than $n$. Then 
$$
\eta^*(\s)\in H^0\big(\X,\,\Om_{\X}^{n-1}\big).
$$
\end{thmn}

\subsection{Reduction steps}
First we would like to recall that, since proving the Extension Theorem is a local problem, it is enough to prove Theorem \ref{t:3} when $V$ is a vector space, $G$ a reductive group acting linearly on $V$ and $X$ is the induced affine GIT-quotient (see Corollary \ref{p:redtoV}).
In this chapter we will reduce the theorem above to an extension problem concerning a single blow up step in the partial resolution of Kirwan. Since for $\dim(X)\leq 1$ the result is clear (in this cases $X$ is smooth), from now on we will assume that $n\geq 2$. To be able to use an induction step we will work in the following setting:

\vspace{10pt}
\fboxsep3mm
\noindent\fbox{
\begin{minipage}[c]{0.945\textwidth}
\begin{setn}
Let $G$ be a reductive group and $V$ a smooth $G$-variety admitting a GIT-quotient $X:=V\hq G$ (with respect to some $G$-linearisation). Assume that for $V$ the partial resolution algorithm of Kirwan works and that we get the following commutative diagram (we know from (K$3$) or (K$3$v) respectively):
\begin{align*}
\xymatrix{
**[l]V^{ss}_k \ar[r]_{F_k^{ss}}\ar@/^1pc/[rrrr]^{F^{ss}}\ar[d]_{\pi_{k}}  & V_{k-1}^{ss} \ar[r]_{F_{k-1}^{ss}}\ar[d]_{\pi_{k-1}} & \dots \ar[r]_{F_{2}^{ss}} & V_1^{ss} \ar[r]_{F_{1}^{ss}}\ar[d]_{\pi_{1}} & **[r] V_0^{ss}:=V^{ss}\ar[d]_{\pi_0} \\
**[l]X_k \ar[r]^{f_k}\ar@/_1pc/[rrrr]_{f}  & X_{k-1} \ar[r]^{f_{k-1}} & \dots \ar[r]^{f_{2}} & X_1 \ar[r]^{f_{1}} & **[r]  X_0:=X.
}
\end{align*}

\end{setn}
\end{minipage}
}
\vspace{10pt}

\begin{lem}
In Setting 1 there exists a log resolution $\eta:\X\rightarrow X$ that factors through the partial resolution of Kirwan.
\end{lem}

\proof
Consider a strong log resolution of $X_k$ and a strong log resolution of $X$ Then there exists a resolution $\eta:\X\rightarrow X$ that dominates both and fulfils the conditions of the lemma. 
\eproof

\vspace{10pt}
\fboxsep3mm
\noindent\fbox{
\begin{minipage}[c]{0.945\textwidth}
\begin{setn}\label{R1}
Let $G$, $V$ and $X:=V\hq G$ be as in Setting 1 and $\eta:\X\rightarrow X$ the resolution introduced in the previous lemma. Then we get the following diagram:
\begin{align}\tag{R1}
\xymatrix{
 & V_k^{ss} \ar[r]_{F_k^{ss}}\ar@/^1pc/[rrrr]^{F^{ss}}\ar[d]_{\pi_{k}}  & V_{k-1}^{ss} \ar[r]_{F_{k-1}^{ss}}\ar[d]_{\pi_{k-1}} & \dots \ar[r]_{F_{2}^{ss}} & V_1^{ss} \ar[r]_{F_{1}^{ss}}\ar[d]_{\pi_{1}} & **[r] V_0^{ss}:=V^{ss}\ar[d]_{\pi_0} \\
 \X\ar[r]^{\eta_k}\ar@/_2pc/[rrrrr]_{\eta} & X_k \ar[r]^{f_k}\ar@/_1pc/[rrrr]_{f}  & X_{k-1} \ar[r]^{f_{k-1}} & \dots \ar[r]^{f_{2}} & X_1 \ar[r]^{f_{1}} & **[r]  X_0:=X.
}
\end{align}

\end{setn}
\end{minipage}
}
\vspace{10pt}

\begin{rem}
The composition $\eta_i:=\eta_k\circ f_k\circ\dots\circ f_{i}:\X\rightarrow X_i$ is a log resolution of $X_i$, for $i\in\{1,\dots,k\}$.
\end{rem}

\reduction{1}
Because of the Corollary \ref{c:Kirw ex} and the invariance of the choice of resolution (see Corollary \ref{c:choice of res}), Theorem \ref{t:3} is a direct corollary of the following proposition:

\begin{prop}\label{p:red 1}
Let $G$, $V$, $X$ and $\eta:\X\rightarrow X$ be as in Setting 2 together with the diagram $(R1)$. Assume that $\dim(X)=n\geq 2$ and that the Extension Theorem is true for all reflexive $p$-forms on any GIT-quotient of dimension less than $n$. Let $\s\in H^0\big(X,\,\Om_{X}^{[n-1]}\big)$ be a reflexive $(n-1)$-form on $X$. Then 
$$
\eta^*(\s)\in H^0\big(\X,\,\Om_{\X}^{n-1}\big).
$$
\end{prop}

\reduction{2}
Recall that $X_k$ has finite quotient singularities. By Proposition \ref{p:ExtThm fqs} the pull-back of any reflexive differential form on $X_k$ by $\eta_k$ extends as a regular form to all of $\X$. Thus $\eta^*(\s)=\eta_k^*(f^*(\s))\in H^0\big(\X,\,\Om_{\X}^{n-1}\big)$ exactly if $f^*(\s)\in H^0\big(X_k,\,\Om_{X_k}^{[n-1]}\big)$. In other words, it is enough to show that the pull-back of $\s$ via $f$ is a reflexive differential form without poles on $X_k$.\\

\reduction{3}\label{r:red 3}
We want to show that in order to prove Proposition \ref{p:red 1} it is enough to prove the following lemma:
\begin{lem}\label{l:red 3}
Let $G$ be a reductive group and $V$ a smooth $G$-variety admitting a GIT-quotient $X:=V\hq G$ (with respect to some $G$-linearisation). Assume that for $V$ the partial resolution algorithm of Kirwan works and that we get the following diagram (that we know from Setting $2$)
\begin{align}\tag{R1}
\xymatrix{
 & V_k^{ss} \ar[r]_{F_k^{ss}}\ar@/^1pc/[rrrr]^{F^{ss}}\ar[d]_{\pi_{k}}  & V_{k-1}^{ss} \ar[r]_{F_{k-1}^{ss}}\ar[d]_{\pi_{k-1}} & \dots \ar[r]_{F_{2}^{ss}} & V_1^{ss} \ar[r]_{F_{1}^{ss}}\ar[d]_{\pi_{1}} & **[r] V_0^{ss}:=V^{ss}\ar[d]_{\pi_0} \\
 \X\ar[r]^{\eta_k}\ar@/_2pc/[rrrrr]_{\eta} & X_k \ar[r]^{f_k}\ar@/_1pc/[rrrr]_{f}  & X_{k-1} \ar[r]^{f_{k-1}} & \dots \ar[r]^{f_{2}} & X_1 \ar[r]^{f_{1}} & **[r]  X_0:=X,
}
\end{align}
where $\eta:\X\rightarrow X$ is a log resolution of $X$, such that the composition $\eta_i:=\eta_k\circ f_k\circ\dots\circ f_{i}:\X\rightarrow X_i$ is a log resolution of $X_i$, for $i\in\{1,\dots,k\}$. Assume that $\dim(X)=n\geq 2$ and that the Extension Theorem is true for all reflexive $p$-forms on any GIT-quotient of dimension less than $n$. Let $\s\in H^0\big(X,\,\Om_{X}^{[n-1]}\big)$ be a reflexive $(n-1)$-form on $X$. Then $\s$ lifts without poles to all of $X_1$ under the first Kirwan step $f_1$. In other words  
$$
\s_1:=f_1^*(\s)\in H^0\big(X_1,\,\Om_{X_1}^{[n-1]}\big).
$$
\end{lem}
Assume that the assertion of Lemma \ref{l:red 3} is true. Then $G$, $V_1$ and $X_1$ fulfil the conditions of Setting $2$ and we get a commutative diagram:

\begin{align}\tag{R2}
\xymatrix{
 & V_k^{ss} \ar[r]_{F_k^{ss}}\ar@/^1pc/[rrr]^{F^{ss}}\ar[d]_{\pi_{k}}  & V_{k-1}^{ss} \ar[r]_{F_{k-1}^{ss}}\ar[d]_{\pi_{k-1}} & \dots \ar[r]_{F_{2}^{ss}} & **[r] V_1^{ss}\ar[d]^{\pi_1:=\pi}\\
 \X\ar[r]^{\eta_k}\ar@/_2pc/[rrrr]_{\eta} & X_k \ar[r]^{f_k}\ar@/_1pc/[rrr]_{f}  & X_{k-1} \ar[r]^{f_{k-1}} & \dots \ar[r]^{f_{2}} & **[r]X_1.
}
\end{align}
By renaming the elements of the diagram we can apply Lemma \ref{l:red 3} and it follows that
$$
f_2^*(\s_1):=f_2^*(f_1^*(\s))\in H^0\big(X_2,\,\Om_{X_2}^{[n-1]}\big).
$$
Continuing with this procedure after finitely many steps we deduce that
$$
f_k^*f_{k-1}^*...f_1^*(\s):=f^*(\s)\in H^0\big(X_k,\,\Om_{X_k}^{[n-1]}\big).
$$
In conclusion, using Reduction $2$ we prove Proposition \ref{p:red 1} by repeatedly applying Lemma \ref{l:red 3} to the GIT-quotients that appear in the partial resolution of Kirwan. 

\subsection{The residue sequence for pairs of good quotients}\label{subs:residuemap}
Let $Y$ be a smooth variety and $D\subset Y$ a reduced irreducible snc divisor. We call $(Y,\, D)$ an \emph{snc pair}. For these pairs \citep[Prop.~2.3]{EV92} introduced the following exact sequence
$$
0\rightarrow \Om^p_{Y}\rightarrow \Om^p_{Y}(\slog D)\xrightarrow{res} \Om^{p-1}_D\rightarrow 0,
$$
for all $1\leq p \leq \dim(Y)$, which is called \emph{residue sequence}. In the case where $Y$ is not smooth the residue sequence was constructed for a dlt pair $(Y,\, D)$ by \citep[Thm.~11.7]{GKKP11}. Assume that $f:Y\rightarrow X$ is a resolution of singularities of $X$ where $D$ is the $f$-exceptional divisor and $\s$ is a reflexive $p$-form on $X$. Assume that $f^*(\s)\in H^0\big(Y,\,\Om_{Y}^{p}(\slog D)\big)$. Then the residue sequence is a useful tool the examine whether $f^*(\s)$ has logarithmic poles along $D$ or not. To this end, consider the induced long exact sequence of cohomology groups
$$
0\rightarrow H^0\big(Y,\,\Om_{Y}^{p}\big)\rightarrow H^0\big(Y,\,\Om_{Y}^{p}(\slog D)\big)\xrightarrow{res} H^0\big(D,\,\Om_{D}^{p-1}\big)\rightarrow \dots~.
$$
Then $f^*(\s)\in H^0\big(Y,\,\Om_{Y}^{p}\big)$ if and only if $res(f^*(\s))=0$ (by abuse of notation we will not differ between the residue map on the level of sheaves or cohomology groups). The following lemma enables us to use a similar approach for GIT-quotients.

\begin{lem}[Residue sequence for good quotients]\label{l:residue map}
Let $G$ be a reductive group and $V$ a smooth $G$-variety admitting a good quotient $X:=V\hq G$. Let $E\subset X$ be a reduced irreducible smooth divisor on $X$. Then there exists a residue sequence (which is exact)
$$
0\rightarrow \Om^{[p]}_X\rightarrow \Om^{[p]}_{X}(\slog E)\xrightarrow{res} \Om^{p-1}_E
$$
for all $1\leq p\leq \dim(X)$, such that the restriction to the sheaves defined on the smooth locus $X_{\mathrm{sm}}$ provides the exact residue sequence of the snc pair $(X_{\mathrm{sm}},E\cap X_{\mathrm{sm}})$:
$$
0\rightarrow \Om^p_{X_{\mathrm{sm}}}\rightarrow \Om^p_{X_{\mathrm{sm}}}(\slog (E\cap X_{\mathrm{sm}}))\xrightarrow{res_{\mid X_{\mathrm{sm}}}} \Om^{p-1}_{E\cap X_{\mathrm{sm}}}\rightarrow 0.
$$
\end{lem}

\proof
Let $i:X_{\mathrm{sm}}\rightarrow X$ be the natural embedding of the smooth locus $X_{\mathrm{sm}}$ into $X$. Consider the exact sequence 
$$
0\rightarrow \Om^p_{X_{\mathrm{sm}}}\rightarrow \Om^p_{X_{\mathrm{sm}}}(\slog (E\cap X_{\mathrm{sm}}))\xrightarrow{res_{\mid X_{\mathrm{sm}}}} \Om^{p-1}_{E\cap X_{\mathrm{sm}}}\rightarrow 0.
$$
We now can use the fact that push-forward is a left-exact functor to obtain the  following exact sequence:
$$
0\rightarrow i_*\Om^p_{X_{\mathrm{sm}}}\rightarrow i_*\Om^p_{X_{\mathrm{sm}}}(\slog (E\cap X_{\mathrm{sm}}))\xrightarrow{res_{\mid X_{\mathrm{sm}}}} i_*\Om^{p-1}_{E\cap X_{\mathrm{sm}}}
$$
By definition of the sheaf of reflexive differential forms, this is the sequence 
$$
0\rightarrow \Om^{[p]}_X\rightarrow \Om^{[p]}_{X}(\slog E)\xrightarrow{res} \Om^{p-1}_E
$$
with the claimed property.
\eproof

\subsection{Log-poles along the exceptional divisor in the first Kirwan step}
Our aim is to prove Lemma \ref{l:red 3}. To this end we analyse the properties of the $f_1$-exceptional divisor $E_1\subset X_1$ that arises in the first Kirwan step. 

\begin{lem}\label{l:pole on birat transf}
Let $X$ be a normal variety and $f: Y\rightarrow X$ a resolution of singularities. Let $E\subset X$ be a reduced divisor on $X$ and $D\subset Y$ the strict transform of $E$ in $Y$. Consider a reflexive differential form $\s$ on $X$. Then $\s$ has no pole (a log-pole, a pole of degree $1$ or bigger) along $E$ if and only if the pull-back $f^*(\s)$ has no pole (a log-pole, a pole of degree $1$ or bigger) along $D$.
\end{lem}

\proof
Since $X$ is normal there exists a point $p\in E\cap X_{\mathrm{sm}}$ and a neighbourhood $U\subset X_{\mathrm{sm}}$ of $p$ such that $f$ is an isomorphism over $U$. The result of the lemma is true for $\s_{\mid U}$ on $E\cap U$ and thus follows for $\s$ on $E$. 
\eproof

\begin{cor}\label{c:log on E1}
In the setting of Lemma \ref{l:red 3}, let $E_1\subset X_1$ be the $f_1$-exceptional divisor. Then 
$$\s_1:=f_1^*(\s)\in H^0\big(X_1,\,\Om_{X_1}^{[n-1]}(\slog E_1)\big).$$
\end{cor}

\proof
Let $\widetilde{E}=\exc(\eta)\subset \X$ be the exceptional divisor on $\X$. By Proposition \ref{logpoles} $\eta^*(\s)\in H^0\big(\X,\,\Om_{\X}^{n-1}(\slog \widetilde{E})\big)$. Recall that $\widetilde E=E_{\eta_k}\cup(\bar E_1\cup \dots \cup \bar E_k)$ where $E_{\eta_k}$ gets contracted by $\eta_k$ and $\bar E_i$ is the irreducible component of $\widetilde E$, that is the strict transform of the $f_i$-exceptional divisor $E_i\subset X_i$, for $i\in\{1,\dots,k\}$. Consider the log resolution $\eta_i:\X\rightarrow X_i$ for $i\in\{1,\dots,k\}$. Since $\eta^*(\s)=\eta_1^*(\s_1)$ has at most a log-pole along $\bar E_1$ it follows that 
$$\s_1:=f_1^*(\s)\in H^0\big(X_1,\,\Om_{X_1}^{[n-1]}(\slog E_1)\big).$$
\eproof

\section{The proof}\label{ch:proofof3}
In this chapter we want to apply the residue sequence constructed in Lemma \ref{l:residue map} to the single Kirwan step (that we consider in Lemma \ref{l:red 3}) to prove Lemma \ref{l:red 3} and thus conclude Theorem \ref{t:3}. As a last preparation we need to introduce the notion of an unirational variety.

\begin{defn}[{\citep[p.~208]{Sha77}}]
A variety $X$ is called \emph{unirational} if there exists a dominant rational morphism $\varPhi:\PP^m\dashrightarrow X$, for some $m\in\NN$.
\end{defn} 

\begin{ex}\label{P/H unirat}
Let $G$ be a reductive group acting on $V=\PP^m$. Consider a $G$-linearisation on $V$ such that we get a projective GIT-quotient $X:=V\hq G$. Since the space of semi-stable points $V^{ss}$ is dense in $V$, the quotient $X$ is a unirational variety.
\end{ex}

\begin{lem}\label{l:unirational}
Let $X$ be a normal unirational variety and let $0<p\leq \dim(X)$ be an integer. Assume that for any reflexive $p$-form $\s\in H^0\big(X,\,\Om^{[p]}_X\big)$ on $X$ and any resolution $f:Y\rightarrow X$ the pull-back $f^*(\s)$ extends to all of $Y$ without poles. Then $H^0\big(X,\,\Om^{[p]}_X\big)=0$.
\end{lem}

\proof
Let $\s \in H^0\big(X,\,\Om^{[p]}_X\big)$ be any reflexive $p$-form on $X$, for $0\leq p\leq dim(X)$, and let $f:Y\rightarrow X$ be a resolution. Then $f^*(\s)\in H^0\big(Y,\,\Om^{p}_{Y}\big)$. Since $X$ is unirational, $Y$ is unirational as well. By \citep[p.~208]{Sha77} we know that $H^0\big(Y,\,\Om^{p}_{Y}\big)=0$ for all $0<p\leq \dim(Y)$. Thus $f^*(\s)=0$ and hence $\s=0$.
\eproof

By Reduction $3$ in Chapter \ref{r:red 3} we know that Theorem \ref{t:3} can be deduced from Lemma \ref{l:red 3}. Our goal is to prove this lemma. For the convenience of the reader we would like to sum up the essentials of the setting we use in the proof:

\vspace{10pt}
\fboxsep3mm
\noindent\fbox{
\begin{minipage}[c]{0.945\textwidth}
\begin{setn}
Let $G$, $V$ and $X:=V\hq G$, with $\dim(X)=:n\geq 2$, be as in Lemma \ref{l:red 3}. Then we get the following diagram
\begin{align}\tag{R1}
\xymatrix{
 & V_k^{ss} \ar[r]_{F_k^{ss}}\ar@/^1pc/[rrrr]^{F^{ss}}\ar[d]_{\pi_{k}}  & V_{k-1}^{ss} \ar[r]_{F_{k-1}^{ss}}\ar[d]_{\pi_{k-1}} & \dots \ar[r]_{F_{2}^{ss}} & V_1^{ss} \ar[r]_{F_{1}^{ss}}\ar[d]_{\pi_{1}} & **[r] V_0^{ss}:=V^{ss}\ar[d]_{\pi_0} \\
 \X\ar[r]^{\eta_k}\ar@/_2pc/[rrrrr]_{\eta} & X_k \ar[r]^{f_k}\ar@/_1pc/[rrrr]_{f}  & X_{k-1} \ar[r]^{f_{k-1}} & \dots \ar[r]^{f_{2}} & X_1 \ar[r]^{f_{1}} & **[r]  X_0:=X,
}
\end{align}
where $\eta:\X\rightarrow X$ is a log-resolution. Let $E_1\subset X_1$ be the $f_1$-exceptional divisor in $X_1$ and let $\s\in H^0\big(X,\,\Om^{[n-1]}_X\big)$ be a reflexive $(n-1)$-form on $X$. Then we want to show that 
$$
\s_1:=f_1^*(\s)\in H^0\big(X_1,\,\Om_{X_1}^{[n-1]}\big).
$$
By Corollary \ref{c:log on E1} we know that $\s_1\in H^0\big(X_1,\,\Om_{X_1}^{[n-1]}(\slog E_1)\big)$. Thus we only have to show that $\s_1$ has no log-pole along $E_1$.

\end{setn}
\end{minipage}
}

\vspace{10pt}

\reduction{4}
Recall that there exists a closed subset $Z\subset X$ with $\codim_X(Z)\geq 3$, such that $X\bsl Z$ has finite quotient singularities (see Proposition \ref{p:codim2} and remark). Set $T:=X_{\mathrm{sing}}\cap Z \subset X$ to be the closed subset of the singular locus that contains all other points. We have to analyse two cases. If $f_1^{-1}(T)\cap E_1\subset X_1$ is not a divisor in $X_1$ then the image of a generic point of $E_1$ lies in $X\bsl T$. In this case we can show that $\s_1$ has no pole along $E_1$, using Corollary \ref{c:codim2} and Lemma \ref{l:pole on birat transf}. 

Hence let us assume for the remainder of the proof that $f_1^{-1}(T)\cap E_1\subset X_1$ is a divisor in $X_1$. Then it is enough to show that $\s_1$ has no pole in a generic point of the divisor $f_1^{-1}(T)\cap E_1$ (see proof of Corollary \ref{c:Extdim3}). Thus, after shrinking $T$, we may assume that $T$ is irreducible of $\codim_X(T)=:l\geq 3$, $E_1=f_1^{-1}(T)$ and that $f_{1\mid E_1}:E_1\rightarrow T$ is a surjective morphism. Since the extension of differential forms is a local problem we may shrink $T$ even further, such that $T$ is smooth and such that the sheaf $\Om_T^q$ of $q$-forms on $T$ ($0\leq q\leq \dim(T)$) is trivial for all $0\leq q\leq \dim(T)$.

We now can prove Lemma \ref{l:red 3} by analysing different dimensions of $T$ separately. Recall that we work under the following assumption.

\vspace{10pt}
\fboxsep3mm
\noindent\fbox{
\begin{minipage}[c]{0.945\textwidth}
\textbf{Assumption (A1).}
The Extension Theorem is true for all reflexive $p$-forms on any GIT-quotient of dimension less than $n$.\label{A1}
\end{minipage}
}
\vspace{10pt}

\sub{Claim 1} 
There exists a closed subspace $Z\subset X_1$ with $\codim_{X_1}(Z)\geq 3$ such that $X_1^0:=X_1\bsl Z$ has finite quotient singularities and such that $E_1^0:=E_1\cap X_1$ is a smooth irreducible divisor on $X_1$.\\

\proofof{Claim 1}
Since $X_1$ and $E_1$ are GIT-quotients of smooth varieties with respect to some $G$-linearisation (see Lemma \ref{l:fiber of Kirwan Step}), they are normal varieties, that are smooth outside codimension $2$ and have only f.q.s. outside codimenion $3$ by Proposition \ref{p:codim2}. \hfill Q.E.D.\\

\noindent\textbf{Case 1 ($\dim(T)=0)$:}
Let $T=\{x\}$. Then $E_1=f_1^{-1}(x)$ is the fiber of $x$ in $X_1$. By Lemma \ref{l:fiber of Kirwan Step} there exists a $m\in \NN_{\geq 1}$ and a reductive subgroup $H\subset G$, such that $E_1$ is the projective GIT-quotient $\PP^m\hq H$ with respect to some $H$-linearisation on $\PP^m$. Thus $E_1$ is unirational by Example \ref{P/H unirat}. \\

Since $E_1$ is unirational and by Lemma \ref{l:fiber of Kirwan Step} a projective GIT-quotient, assumption (A$1$) ensures that we can use Lemma \ref{l:unirational} to get the following result:
\begin{align}\tag{+}
H^0\big(E_1,\,\Om^{[p]}_{E_1}\big)=0 \text{ for } 0<p\leq \dim(E_1)=n-1.
\end{align}

Let $Z$, $X_1^0$, and $E_1^0$ be as in Claim $1$. Recall that $\s_1\in H^0\big(X_1,\,\Om_{X_1}^{[n-1]}(\slog E_1)\big)$. Thus the restriction $\s_1^0:=\s_{1\mid X_1^0}\in H^0\big(X_1^0,\,\Om_{X_1^0}^{[n-1]}(\slog E_1^0)\big)$ might have a log pole along $E_1^0$. We can use the residue sequence constructed in Lemma \ref{l:residue map} to define a $(n-1)$-form $\alpha:=res(\s_1^0)$ on $E_1^0$. Since $E_1^0=E_1\bsl Z'$, for a closed subset $Z'=Z\cap E_1\subset E_1$ of $\codim_{E_1}(Z')\geq 2$ and $E_1$ is normal, we see that (+) implies that $H^0\big(E_1^0,\,\Om^{n-2}_{E_1^0}\big)=0$. Thus $\alpha=0$ and by Lemma \ref{l:residue map} hence $\s_1^0\in H^0\big(X_1^0,\,\Om_{X_1^0}^{[n-1]}\big)$. Since $X_1^0=X_1\bsl Z$ with $\codim_{X_1}(Z)\geq 3$ and $X_1$ is normal, it follows that $\s_1\in H^0\big(X_1,\,\Om_{X_1}^{[n-1]}\big)$. This ends the proof of Case $1$. \\

\noindent\textbf{Case 2 ($\dim(T)>0)$:}
Let $t\in T$ be a general point of $T$ and the fiber $f_1^{-1}(t)=:E_{1,t}\subset E_1$. By Lemma \ref{l:fiber of Kirwan Step} there exists a $m_t\in \NN_{\geq 1}$ and a reductive subgroup $H_t\subset G$, such that $E_{1,t}$ is the projective GIT-quotient $\PP^{m_t}\hq H_t$ (with respect to some $H_t$-linearisation on $\PP^{m_t}$). Thus $E_{1,t}$ is unirational. \\

Since for all $t\in T$ the fiber $E_{1,t}$ is unirational and by Lemma \ref{l:fiber of Kirwan Step} a projective GIT-quotient (and therefore normal), assumption (A$1$) ensures that we can use Lemma \ref{l:unirational} to see that $H^0\big(E_{1,t},\,\Om^{[p]}_{E_{1,t}}\big)=0$ for $0<p\leq \dim(E_{1,t})$ and for all $t\in T$. Thus, by the same argument as in \citep[Lem.~25.11]{GKKP11}, we get the following result:
\begin{align}\tag{++}
H^0\big(E_1,\,\Om^{[p]}_{E_1/T}\big)=0 \text{ for }0<p\leq dim (E_1).
\end{align}

Let $Z$, $X_1^0$ and $E_1^0$ be as in Claim $1$. As before $E_1^0=E_1\bsl Z'$, for a closed subset $Z'=Z\cap E_1\subset E_1$ of $\codim_{E_1}(Z')\geq 2$. Then $H^0\big(E_1,\,\Om^{[p]}_{E_1/T}\big)=0$ implies that $H^0\big(E_1^0,\,\Om^{p}_{E_1^0/T}\big)=0$.

Assume that $\s_1\in H^0\big(X_1,\,\Om_{X_1}^{[n-1]}(\slog E_1)\big)$ has a true log pole along $E_1$. Then the restriction $\s_1^0:=\s_{1\mid X_1^0}\in H^0\big(X_1^0,\,\Om_{X_1^0}^{[n-1]}(\slog E_1^0)\big)$ has a log pole along $E_1^0$. We can use the residue sequence constructed in Lemma \ref{l:residue map} to define a non-vanishing $(n-2)$-form $\alpha:=res(\s_1^0)\neq 0$ on $E_1^0$. We will now show that such a form cannot exists.

Consider the induced morphism $\varPsi: E_1^0\rightarrow T$ of smooth varieties. After shrinking $T$, we may assume that $\varPsi$ is a smooth morphism. Thus, we get the following filtration (see \cite[II.~Ex.~5.16]{Har77})
$$
{\Om}^{n-2}_{E_1^0}=\mathcal{F}^0\supset\mathcal{F}^1\supset\dots\supset\mathcal{F}^{n-2}\supset\mathcal{F}^{n-1}=\{0\}.
$$
These filtrations induce exact sequences
\begin{align}\tag{$S_r$}
0\rightarrow \mathcal{F}^{r+1}\rightarrow \mathcal{F}^r \rightarrow \varPsi^*\Om_{T}^{r}\otimes \Om^{n-2-r}_{E_1^0/T}\rightarrow 0,
\end{align}
for all $0\leq r\leq n-2$, where, by construction, $\mathcal{F}^{n-2}\cong \varPsi^*\Om_{T}^{n-2}$. If we consider the long exact sequence induced by the sequence above, the following is true for $0\leq r\leq n-2$: Assume that $\beta$ is a non-vanishing section of $\mathcal F^r$. Then either it induces a non-vanishing section of $\varPsi^*\Om_{T}^{r}\otimes \Om^{n-2-r}_{E_1^0/T}$, or it comes from a non-vanishing section on $\mathcal{F}^{r+1}$. Chasing the $(n-2)$-form $\alpha$ through these sequences we get the following contradiction:

Recall that $\varPsi^*\Om_{T}^{q}$ is a trivial bundle for all $0\leq q\leq \dim(T)$. Thus, by (++) we see, that $H^0\big(E_1^0,\,\varPsi^*\Om_{T}^{r}\otimes \Om^{n-2-r}_{E_1^0/T}\big)=0$ for all $0<n-2-r\leq  \dim(E_1)$. As a conclusion the last sheaf in the sequence ($S_r$) has no non-vanishing section for $0\leq r < n-2$. Therefore, the non-vanishing $(n-2)$-form $\alpha$ on $E_1^0$ induces a non-vanishing 
$(n-2)$-form $\alpha'\in H^0\big(E_1^0,\,\varPsi^*\Om_{T}^{n-2}\big)$. \label{induction Sr} However, $\dim(T)\leq\dim(X)-3=n-3$ and such a form cannot exist. As a consequence $\s_1$ has no log pole along $E_1$ and $\s_1\in H^0\big(X_1,\,\Om_{X_1}^{[n-1]}\big)$.
\eproof

\newpage

\thispagestyle{plain}

\bibliography{literaturs}

\newcommand{\etalchar}[1]{$^{#1}$}
\begin{thebibliography}{GKKP11}


\providecommand{\url}[1]{\texttt{#1}}
\expandafter\ifx\csname urlstyle\endcsname\relax
  \providecommand{\doi}[1]{doi: #1}\else
  \providecommand{\doi}{doi: \begingroup \urlstyle{rm}\Url}\fi

\bibitem[BB{\'S}97]{BS91}
\textsc{Bia{\l}ynicki-Birula}, Andrzej ; \textsc{{\'S}wi{\k e}cicka}, Joanna:
\newblock Three theorems on existence of good quotients.
\newblock {In: }\emph{Math. Ann.} 307 (1997), Nr. 1, 143--149.
\newblock \url{http://dx.doi.org/10.1007/s002080050027}. --
\newblock DOI 10.1007/s002080050027. --
\newblock ISSN 0025--5831

\bibitem[Bou87]{Bou87}
\textsc{Boutot}, Jean-Fran\c{c}ois:
\newblock Singularit\'es rationnelles et quotients par les groupes r\'eductifs.
\newblock {In: }\emph{Invent. Math.} 88 (1987), Nr. 1, 65--68.
\newblock \url{https://doi.org/10.1007/BF01405091}. --
\newblock ISSN 0020--9910

\bibitem[CLS11]{CLS11}
\textsc{Cox}, David~A. ; \textsc{Little}, John~B.  ; \textsc{Schenck},
  Henry~K.:
\newblock \emph{Graduate Studies in Mathematics}. Bd. 124: {\emph{Toric
  varieties}}.
\newblock American Mathematical Society, Providence, RI, 2011. --
\newblock  xxiv+841 S.
\newblock \url{http://dx.doi.org/10.1090/gsm/124}.
\newblock \url{http://dx.doi.org/10.1090/gsm/124}. --
\newblock ISBN 978--0--8218--4819--7

\bibitem[Dol03]{Dol03}
\textsc{Dolgachev}, Igor:
\newblock \emph{London Mathematical Society Lecture Note Series}. Bd. 296:
  {\emph{Lectures on invariant theory}}.
\newblock Cambridge University Press, Cambridge, 2003. --
\newblock  xvi+220 S.
\newblock \url{http://dx.doi.org/10.1017/CBO9780511615436}.
\newblock \url{http://dx.doi.org/10.1017/CBO9780511615436}. --
\newblock ISBN 0--521--52548--9

\bibitem[Dr{\'e}04]{Dre12}
\textsc{Dr{\'e}zet}, Jean-Marc:
\newblock Luna's slice theorem and applications.
\newblock {In: }\emph{Algebraic group actions and quotients}.
\newblock Hindawi Publ. Corp., Cairo, 2004, S. 39--89

\bibitem[EV92]{EV92}
\textsc{Esnault}, H\'el\`ene ; \textsc{Viehweg}, Eckart:
\newblock \emph{DMV Seminar}. Bd.~20: {\emph{Lectures on vanishing theorems}}.
\newblock Birkh\"auser Verlag, Basel, 1992. --
\newblock  vi+164 S.
\newblock \url{http://dx.doi.org/10.1007/978-3-0348-8600-0}.
\newblock \url{http://dx.doi.org/10.1007/978-3-0348-8600-0}. --
\newblock ISBN 3--7643--2822--3

\bibitem[FGI{\etalchar{+}}05]{FGI05}
\textsc{Fantechi}, Barbara ; \textsc{G\"ottsche}, Lothar ; \textsc{Illusie},
  Luc ; \textsc{Kleiman}, Steven~L. ; \textsc{Nitsure}, Nitin  ;
  \textsc{Vistoli}, Angelo:
\newblock \emph{Mathematical Surveys and Monographs}. Bd. 123:
  {\emph{Fundamental algebraic geometry}}.
\newblock American Mathematical Society, Providence, RI, 2005. --
\newblock  x+339 S. --
\newblock ISBN 0--8218--3541--6. --
\newblock Grothendieck's FGA explained

\bibitem[GK14]{GK14}
\textsc{Graf}, Patrick ; \textsc{Kov\'acs}, S\'andor~J.:
\newblock Potentially {D}u {B}ois spaces.
\newblock {In: }\emph{J. Singul.} 8 (2014), S. 117--134. --
\newblock ISSN 1949--2006

\bibitem[GKK10]{GKK10}
\textsc{Greb}, Daniel ; \textsc{Kebekus}, Stefan  ; \textsc{Kov\'acs},
  S\'andor~J.:
\newblock Extension theorems for differential forms and {B}ogomolov-{S}ommese
  vanishing on log canonical varieties.
\newblock {In: }\emph{Compos. Math.} 146 (2010), Nr. 1, 193--219.
\newblock \url{http://dx.doi.org/10.1112/S0010437X09004321}. --
\newblock DOI 10.1112/S0010437X09004321. --
\newblock ISSN 0010--437X

\bibitem[GKKP11]{GKKP11}
\textsc{Greb}, Daniel ; \textsc{Kebekus}, Stefan ; \textsc{Kov\'acs},
  S\'andor~J.  ; \textsc{Peternell}, Thomas:
\newblock Differential forms on log canonical spaces.
\newblock {In: }\emph{Publ. Math. Inst. Hautes \'Etudes Sci.}  (2011), Nr. 114,
  87--169.
\newblock \url{http://dx.doi.org/10.1007/s10240-011-0036-0}. --
\newblock DOI 10.1007/s10240--011--0036--0. --
\newblock ISSN 0073--8301

\bibitem[Gur91]{Gur}
\textsc{Gurjar}, Rajendra~V.:
\newblock On a conjecture of {C}. {T}. {C}. {W}all.
\newblock {In: }\emph{J. Math. Kyoto Univ.} 31 (1991), Nr. 4, 1121--1124.
\newblock \url{http://dx.doi.org/10.1215/kjm/1250519680}. --
\newblock DOI 10.1215/kjm/1250519680. --
\newblock ISSN 0023--608X

\bibitem[GW10]{GW10}
\textsc{G\"ortz}, Ulrich ; \textsc{Wedhorn}, Torsten:
\newblock \emph{Algebraic geometry {I}}.
\newblock Vieweg + Teubner, Wiesbaden, 2010 (Advanced Lectures in Mathematics).
  --
\newblock  viii+615 S.
\newblock \url{http://dx.doi.org/10.1007/978-3-8348-9722-0}.
\newblock \url{http://dx.doi.org/10.1007/978-3-8348-9722-0}. --
\newblock ISBN 978--3--8348--0676--5. --
\newblock Schemes with examples and exercises

\bibitem[Har77]{Har77}
\textsc{Hartshorne}, Robin:
\newblock \emph{Algebraic geometry}.
\newblock Springer-Verlag, New York-Heidelberg, 1977. --
\newblock  xvi+496 S. --
\newblock ISBN 0--387--90244--9. --
\newblock Graduate Texts in Mathematics, No. 52

\bibitem[Har80]{Har80}
\textsc{Hartshorne}, Robin:
\newblock Stable reflexive sheaves.
\newblock {In: }\emph{Math. Ann.} 254 (1980), Nr. 2, 121--176.
\newblock \url{https://doi.org/10.1007/BF01467074}. --
\newblock ISSN 0025--5831

\bibitem[Hau04]{Hau03}
\textsc{Hausen}, J\"urgen:
\newblock Geometric invariant theory based on {W}eil divisors.
\newblock {In: }\emph{Compos. Math.} 140 (2004), Nr. 6, 1518--1536.
\newblock \url{https://doi.org/10.1112/S0010437X04000867}. --
\newblock ISSN 0010--437X

\bibitem[Iit82]{Iit82}
\textsc{Iitaka}, Shigeru:
\newblock \emph{Graduate Texts in Mathematics}. Bd.~76: {\emph{Algebraic
  geometry}}.
\newblock Springer-Verlag, New York-Berlin, 1982. --
\newblock  x+357 S. --
\newblock ISBN 0--387--90546--4. --
\newblock An introduction to birational geometry of algebraic varieties,
  North-Holland Mathematical Library, 24

\bibitem[Keb13]{Keb13}
\textsc{Kebekus}, Stefan:
\newblock Pull-back morphisms for reflexive differential forms.
\newblock {In: }\emph{Adv. Math.} 245 (2013), 78--112.
\newblock \url{http://dx.doi.org/10.1016/j.aim.2013.06.013}. --
\newblock DOI 10.1016/j.aim.2013.06.013. --
\newblock ISSN 0001--8708

\bibitem[Kem77]{Kem77}
\textsc{Kempf}, George~R.:
\newblock Some quotient varieties have rational singularities.
\newblock {In: }\emph{Michigan Math. J.} 24 (1977), Nr. 3, 347--352.
\newblock \url{http://projecteuclid.org/euclid.mmj/1029001952}. --
\newblock ISSN 0026--2285

\bibitem[Kir85]{Kir85}
\textsc{Kirwan}, Frances~C.:
\newblock Partial desingularisations of quotients of nonsingular varieties and
  their {B}etti numbers.
\newblock {In: }\emph{Ann. of Math. (2)} 122 (1985), Nr. 1, 41--85.
\newblock \url{http://dx.doi.org/10.2307/1971369}. --
\newblock DOI 10.2307/1971369. --
\newblock ISSN 0003--486X

\bibitem[KM98]{KM98}
\textsc{Koll\'ar}, J\'anos ; \textsc{Mori}, Shigefumi:
\newblock \emph{Cambridge Tracts in Mathematics}. Bd. 134: {\emph{Birational
  geometry of algebraic varieties}}.
\newblock Cambridge University Press, Cambridge, 1998. --
\newblock  viii+254 S.
\newblock \url{http://dx.doi.org/10.1017/CBO9780511662560}.
\newblock \url{http://dx.doi.org/10.1017/CBO9780511662560}. --
\newblock ISBN 0--521--63277--3. --
\newblock With the collaboration of C. H. Clemens and A. Corti, Translated from
  the 1998 Japanese original

\bibitem[Kov99]{Kov99}
\textsc{Kov\'acs}, S\'andor~J.:
\newblock Rational, log canonical, {D}u {B}ois singularities: on the
  conjectures of {K}oll\'ar and {S}teenbrink.
\newblock {In: }\emph{Compositio Math.} 118 (1999), Nr. 2, 123--133.
\newblock \url{http://dx.doi.org/10.1023/A:1001120909269}. --
\newblock DOI 10.1023/A:1001120909269. --
\newblock ISSN 0010--437X

\bibitem[Kov00]{Kov00}
\textsc{Kov\'acs}, S\'andor~J.:
\newblock A characterization of rational singularities.
\newblock {In: }\emph{Duke Math. J.} 102 (2000), Nr. 2, 187--191.
\newblock \url{http://dx.doi.org/10.1215/S0012-7094-00-10221-9}. --
\newblock DOI 10.1215/S0012--7094--00--10221--9. --
\newblock ISSN 0012--7094

\bibitem[Kra84]{Kra84}
\textsc{Kraft}, Hanspeter:
\newblock \emph{Geometrische {M}ethoden in der {I}nvariantentheorie}.
\newblock Friedr. Vieweg \& Sohn, Braunschweig, 1984 (Aspects of Mathematics,
  D1). --
\newblock  x+308 S.
\newblock \url{http://dx.doi.org/10.1007/978-3-322-83813-1}.
\newblock \url{http://dx.doi.org/10.1007/978-3-322-83813-1}. --
\newblock ISBN 3--528--08525--8

\bibitem[LR79]{LR79}
\textsc{Luna}, Domingo ; \textsc{Richardson}, Roger~W.:
\newblock A generalization of the {C}hevalley restriction theorem.
\newblock {In: }\emph{Duke Math. J.} 46 (1979), Nr. 3, 487--496.
\newblock \url{http://projecteuclid.org/euclid.dmj/1077313569}. --
\newblock ISSN 0012--7094

\bibitem[Lun73]{Lun73}
\textsc{Luna}, Domingo:
\newblock Slices \'etales.
\newblock   (1973), S. 81--105. Bull. Soc. Math. France, Paris, M\'emoire 33

\bibitem[Mat60]{Mat60}
\textsc{Matsushima}, Yoz\^o:
\newblock Espaces homog\`enes de {S}tein des groupes de {L}ie complexes.
\newblock {In: }\emph{Nagoya Math. J} 16 (1960), 205--218.
\newblock \url{http://projecteuclid.org/euclid.nmj/1118800370}. --
\newblock ISSN 0027--7630

\bibitem[MFK94]{MFK94}
\textsc{Mumford}, David ; \textsc{Fogarty}, John  ; \textsc{Kirwan},
  Frances~C.:
\newblock \emph{Ergebnisse der Mathematik und ihrer Grenzgebiete (2) [Results
  in Mathematics and Related Areas (2)]}. Bd.~34: {\emph{Geometric invariant
  theory}}.
\newblock Third.
\newblock Springer-Verlag, Berlin, 1994. --
\newblock  xiv+292 S.
\newblock \url{http://dx.doi.org/10.1007/978-3-642-57916-5}.
\newblock \url{http://dx.doi.org/10.1007/978-3-642-57916-5}. --
\newblock ISBN 3--540--56963--4

\bibitem[Nam01a]{Nam01a}
\textsc{Namikawa}, Yoshinori:
\newblock Deformation theory of singular symplectic {$n$}-folds.
\newblock {In: }\emph{Math. Ann.} 319 (2001), Nr. 3, 597--623.
\newblock \url{http://dx.doi.org/10.1007/PL00004451}. --
\newblock DOI 10.1007/PL00004451. --
\newblock ISSN 0025--5831

\bibitem[Nam01b]{Nam01}
\textsc{Namikawa}, Yoshinori:
\newblock Extension of 2-forms and symplectic varieties.
\newblock {In: }\emph{J. Reine Angew. Math.} 539 (2001), 123--147.
\newblock \url{http://dx.doi.org/10.1515/crll.2001.070}. --
\newblock DOI 10.1515/crll.2001.070. --
\newblock ISSN 0075--4102

\bibitem[Pin77]{Pin77}
\textsc{Pinkham}, Henry:
\newblock Normal surface singularities with {$C\sp*$} action.
\newblock {In: }\emph{Math. Ann.} 227 (1977), Nr. 2, 183--193.
\newblock \url{http://dx.doi.org/10.1007/BF01350195}. --
\newblock DOI 10.1007/BF01350195. --
\newblock ISSN 0025--5831

\bibitem[Ses72]{Ses72}
\textsc{Seshadri}, Conjeervaram~S.:
\newblock Quotient spaces modulo reductive algebraic groups.
\newblock {In: }\emph{Ann. of Math. (2)} 95 (1972), 511--556; errata, ibid. (2)
  96 (1972), 599.
\newblock \url{http://dx.doi.org/10.2307/1970870}. --
\newblock DOI 10.2307/1970870. --
\newblock ISSN 0003--486X

\bibitem[Sha77]{Sha77}
\textsc{Shafarevich}, Igor~R.:
\newblock \emph{Basic algebraic geometry}.
\newblock Study.
\newblock Springer-Verlag, Berlin-New York, 1977. --
\newblock  xv+439 S. --
\newblock Translated from the Russian by K. A. Hirsch, Revised printing of
  Grundlehren der mathematischen Wissenschaften, Vol. 213, 1974

\bibitem[Sha13]{Sha77b}
\textsc{Shafarevich}, Igor~R.:
\newblock \emph{Basic algebraic geometry. 2}.
\newblock Third.
\newblock Springer, Heidelberg, 2013. --
\newblock  xiv+262 S. --
\newblock ISBN 978--3--642--38009--9; 978--3--642--38010--5. --
\newblock Schemes and complex manifolds, Translated from the 2007 third Russian
  edition by Miles Reid

\bibitem[SS85]{vSS85}
\textsc{Straten}, Duco van ; \textsc{Steenbrink}, Joseph:
\newblock Extendability of holomorphic differential forms near isolated
  hypersurface singularities.
\newblock {In: }\emph{Abh. Math. Sem. Univ. Hamburg} 55 (1985), 97--110.
\newblock \url{http://dx.doi.org/10.1007/BF02941491}. --
\newblock DOI 10.1007/BF02941491. --
\newblock ISSN 0025--5858

\bibitem[Ste83]{Ste83}
\textsc{Steenbrink}, Joseph H.~M.:
\newblock Mixed {H}odge structures associated with isolated singularities.
\newblock {In: }\emph{Singularities, {P}art 2 ({A}rcata, {C}alif., 1981)}
  Bd.~40.
\newblock Amer. Math. Soc., Providence, RI, 1983, S. 513--536

\end{thebibliography}

\end{document}